\documentclass[preprint, 11pt]{elsarticle}

\newcommand{\bsub}{\begin{subequations}}
\newcommand{\esub}{\end{subequations}$\!$}
\newcommand{\ds}[0]{\displaystyle}

\usepackage{amsfonts}
\usepackage[fleqn,reqno]{amsmath}
\usepackage{amssymb}
\usepackage[titletoc]{appendix}
\usepackage{array}
\usepackage{enumitem}
\usepackage{filecontents}
\usepackage[top=1.2in,bottom=1.2in,left=1in, right=1in]{geometry}
\usepackage{graphics}
\usepackage{lineno}
\usepackage{pgfplots}
\usepackage{tikz}
\usepackage{todonotes}
\usetikzlibrary{arrows}
\usepackage{comment}

\usepackage[pagebackref=false,bookmarks=false]{hyperref} 

\hypersetup{
  bookmarksnumbered=true,
  bookmarksopen=false,
  hypertexnames=false,      
  breaklinks=true,          
  unicode=false,
  pdffitwindow=true,        
  pdfnewwindow=true,        
  colorlinks=true,         
  linkcolor=dblue,
  anchorcolor=red,
  citecolor=dorange,
  filecolor=magenta,
  urlcolor=dblue,
  pdfstartview = FitH,
  pdfkeywords = {},
  pdfcreator = {LaTeX with hyperref package}
}

\newcommand{\bd}{{\partial}}

\newcommand{\CC}{{\mathbb{C}}}
\newcommand{\DD}{{\mathcal{D}}}

\newcommand{\LL}{{\mathcal{L}}}

\newcommand{\nn}{{\mathbf{n}}}

\newcommand{\RR}{{\mathbb{R}}}

\renewcommand{\SS}{{\mathcal{S}}}

\newcommand{\xx}{{\mathbf{x}}}

\newcommand{\yy}{{\mathbf{y}}}

\newcommand{\pderiv}[2]{\frac{\partial #1}{\partial #2}}

\begin{document}

\title{Boundary Integral Methods for Particle Diffusion in Complex
Geometries: Shielding, Confinement, and Escape}

\author[Bryan1]{Jesse Cherry}
\author[Alan]{Alan E. Lindsay}
\author[Bryan1,Bryan2]{Bryan D.~Quaife}

\address[Bryan1]{Department of Scientific Computing, Florida State University, Tallahassee, FL, 32306.}
\address[Bryan2]{Geophysical Fluid Dynamics Institute, Florida State
University, Tallahassee, FL, 32306.}
\address[Alan]{Department of Applied \& Computational Math \&
Statistics, University of Notre Dame, Notre Dame, IN, 46656.}

\begin{abstract}
We present a numerical method for the solution of diffusion problems in
unbounded planar regions with complex geometries of absorbing and reflecting bodies. Our numerical method applies
the Laplace transform to the parabolic problem, yielding a modified
Helmholtz equation which is solved with a boundary integral method.
Returning to the time domain is achieved by quadrature of the inverse
Laplace transform by deforming along the so-called Talbot contour. We
demonstrate the method for various complex geometries formed by disjoint
bodies of arbitrary shape on which either uniform Dirichlet or Neumann
boundary conditions are applied. The use of the Laplace transform
bypasses constraints with traditional time-stepping methods and allows
for integration over the long equilibration timescales present in
diffusion problems in unbounded domains. Using this method, we
demonstrate shielding effects where the complex geometry
modulates the dynamics of capture to absorbing sets. In particular, we
show examples where geometry can guide diffusion processes to particular
absorbing sites, obscure absorbing sites from diffusing particles, and
even find the exits of confining geometries, such as mazes. 
\end{abstract}


\maketitle

\section{Introduction\label{s:intro}}

Many problems in biological transport involve a diffusing particle from
a source to a fixed or mobile target~\cite{nic-mul2023,
ste-zmu-hen-rij-ors-lin2020, che-war2024,Lindsay2023a,Lindsay2023b,BN2013,Isaacson2020}. The
target, denoted by $\Gamma_D$, may be a preferred ecological habitat, a
receptor on the surface of a cell, or an appropriate mate. Additionally,
the particle may have to navigate through complex geometries and reflecting obstacles, denoted
by $\Gamma_N$. In the scenario of unbiased diffusion in a planar region
$\Omega \subset \mathbb{R}^2$, this requires solving the stochastic
differential equation
\bsub\label{eqn:SDE}
\begin{equation}\label{eqn:SDE_a}
d \xx = \sqrt{2D} \, d\textbf{W}, \quad t>0, \quad \xx \in\Omega; \qquad
\xx (0) = \xx^{*}\in\Omega.
\end{equation}
Here $\xx^{*}$ is the initial particle location, $D>0$ is the
diffusivity coefficient, and $d\textbf{W}$ is the increment of a Weiner
process. The main quantity of interest is the first arrival time $\tau$
at which the particle arrives at a target region $\Gamma_D$, given by
\begin{equation}\label{eqn:SDE_b}
  \tau = \min_{t>0} \{t \:|\: \xx(t)\in\Gamma_D \}.
\end{equation}
\esub
Direct simulation of~\eqref{eqn:SDE} offers a relatively straightforward
method for sampling the distribution of $\tau$, however, achieving high
accuracy is tricky on account of errors near boundaries, the need to
integrate for long times, and slow convergence of Monte Carlo methods.
These issues can be mitigated with efficiency strategies such as Kinetic
Monte Carlo~\cite{LBS2018,lindsay2024} or Walk on
Spheres~\cite{Northrup1988,HWANG20101089,OKMascagni2004}, however, it is
desirable to have an deterministic solver for key quantities such
as~\eqref{eqn:SDE_b}.
 
\begin{figure}[htbp]
\centering
\includegraphics[width = 0.45\textwidth]{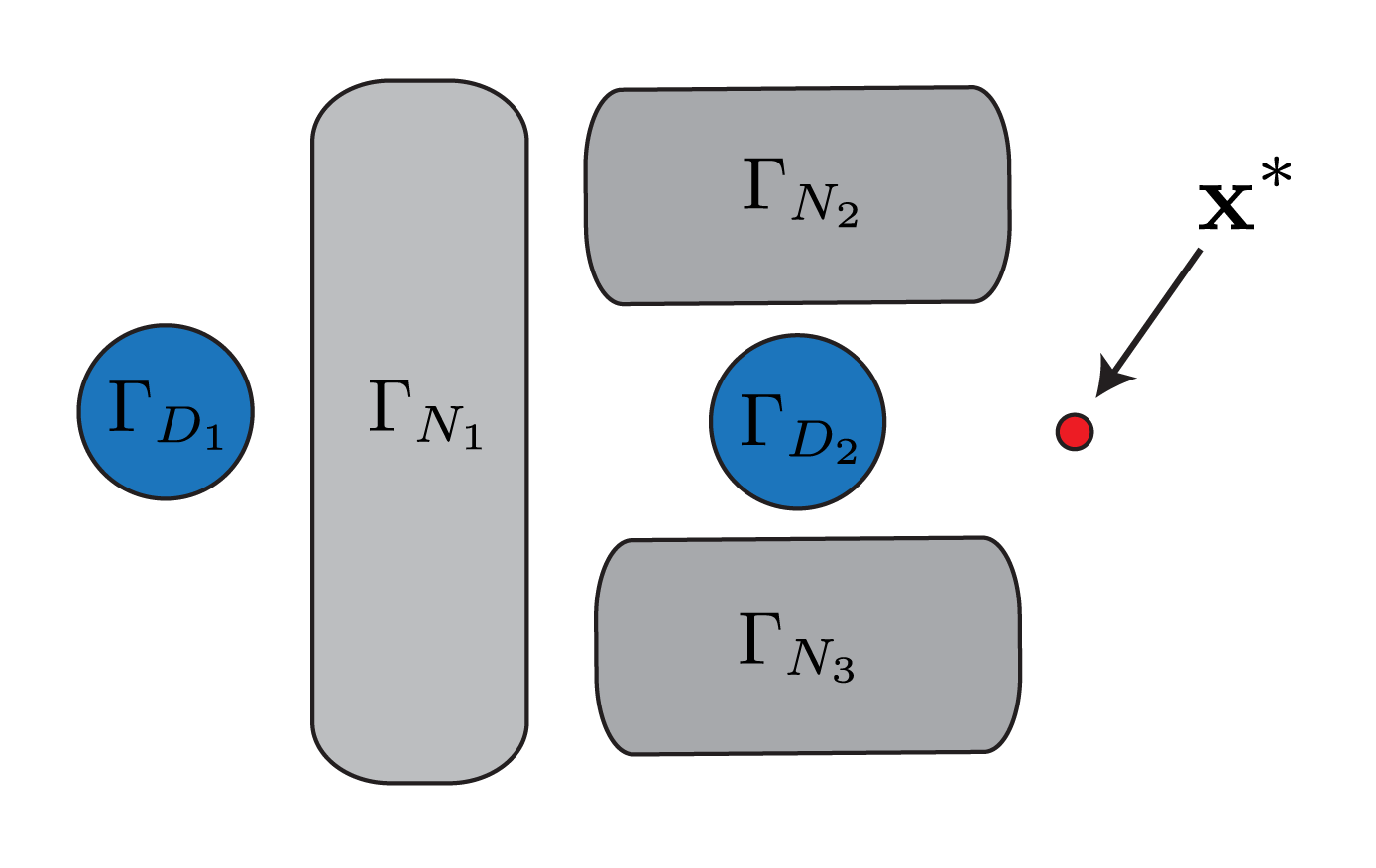}
\caption{A schematic of the domain configuration comprising Neumann
$\Gamma_N = \cup_{k=1}^{M_N} \Gamma_{N_k}$ and Dirichlet $\Gamma_D =
\cup_{k=1}^{M_D} \Gamma_{D_k}$ components and source location
$\xx^{\ast}\in\Omega$. \label{fig:intro}}
\end{figure}

In this work we propose to bypass these issues by solving directly for
the density probability distribution of the particle. The forward
Kolmogorov, or Fokker-Planck equation, describes the
process~\eqref{eqn:SDE} as a parabolic partial differential equation
(PDE). We consider an unbounded multiply-connected domain $\Omega
\subset \RR^2$ with boundary $\Gamma = \bd\Omega$. For simplicity, we
set the diffusivity to be $D=1$ as shown in the schematic of
Figure~\ref{fig:intro}. The boundary condition on each component of
$\Gamma$ is homogeneous and can be either Dirichlet (absorbing bodies)
or Neumann (reflecting bodies). The $M_D$ components with Dirichlet
boundary conditions are denoted by $\Gamma_D = \cup_{k=1}^{M_D}
\Gamma_{D_k}$, and the $M_N$ components with Neumann boundary conditions
are denoted by $\Gamma_N = \cup_{k=1}^{M_N} \Gamma_{N_k}$. Therefore,
$\Gamma = \Gamma_D \cup \Gamma_N$. Note that we do not allow for a
single connected component of $\Gamma$ to be divided into both Dirichlet
and Neumann regions. Then, the governing equation is
\bsub\label{eqn:diffusion}
  \begin{alignat}{3}
\label{eqn:diffusion_a}    \pderiv{p}{t} &= \Delta p, &&\xx \in \Omega,  \\
\label{eqn:diffusion_b}    p &= 0, &&\xx \in \Gamma_D, \\
\label{eqn:diffusion_c}    \pderiv{p}{\nn} &= 0, &&\xx \in \Gamma_N, \\
\label{eqn:diffusion_d}   p(\xx,0) &= \delta(\xx - \xx^*), \qquad &&\xx \in \Omega,
  \end{alignat}
\esub
where $\nn$ is the unit outward normal of $\Gamma$. The initial
condition~\eqref{eqn:diffusion_d} corresponds to the particle initially
at $\xx^* \in \Omega$. We also require that $p$ is bounded as
$|\xx|\rightarrow \infty$. 

A central quantity of interest from the
solution of \eqref{eqn:diffusion} is the flux through the target region,
and its corresponding cumulative flux, defined as
\begin{equation}
  \label{eqn:bdflux}
  j(t) = \int_{\Gamma_D} \pderiv{p}{\nn} \, ds, \qquad
  c(t) = \int_{0}^{t} j(\eta) \, d\eta,
\end{equation}
respectively. The flux $j(t)$ gives the probability distribution of the
first arrival time \eqref{eqn:SDE_b} and typically has a large amount of
mass at late time, so that the cumulative flux $c(t)$ converges to one
very slowly as $t\to\infty$. In terms of the capture time variable $\tau$ defined by \eqref{eqn:SDE}, we have that
\begin{equation}
j(t) = \mathbb{P}(\tau=t), \qquad c(t) = \mathbb{P}(\tau<t).
\end{equation}

In this paper we develop a BIE method to solve the
PDE~\eqref{eqn:diffusion}, and to compute the flux $j(t)$ and cumulative flux $c(t)$ defined by~\eqref{eqn:bdflux}.
Integral equation methods offer a powerful technique over conventional
stencil-based methods for solving PDEs, including: unbounded complex
domains are easily incorporated, appropriately chosen formulations
result in linear systems whose conditioning is mesh-independent, and
linear complexity is possible by incorporating fast summation methods or
fast direct solvers.

One approach to solve equation~\eqref{eqn:diffusion} is to first
discretize the PDE in time, treating the diffusion term implicitly to
reduce numerical stiffness. This results in an elliptic PDE that must be
solved at each time step, and BIEs have been used to solve this PDE by
many groups (Ex.~\cite{cha-kre1997, fry-kro-tor2020, kro-qua2010}). A
similar approach has been applied to other time-dependent PDEs such as
the Navier-Stokes equations~\cite{kli-ask-kro2020} and the wave
equation~\cite{cau-chr-ong-van2014}. Planar diffusion problems are
characterized by long timescale dynamics, hence standard time stepping
approaches are not practical as an excessively small time step size is
required for temporal resolution.

Alternatively, the heat equation can be written in terms of the
time-dependent heat kernel and this results in a Volterra integral
equation. The regularity of the solution, and the compactness, and
coercivity of the integral equations have been analyzed by several
groups (Ex.~\cite{hsi-sar1993, kre1999, qiu-rie-say-zha2017}). A naive
numerical approach to solve the integral equation requires storing and
integrating over the entire time history, but this is computationally
prohibitive. Greengard and Strain~\cite{gre-str1990} introduced an
optimal complexity algorithm that makes use of the smoothness of
the history part of the heat kernel. A necessary algorithm of this
approach is the Gauss transform, and this has led to several works on
computing fast Gauss transforms~\cite{str1994, vee-bir2008,
wan-gre2018}. Several groups perform a change of variables including a
Fourier transform in space or time~\cite{li-gre2007, li-gre2009}, an
Abel transform~\cite{mel-res2014}, a Chebyshev transform~\cite{tau2007,
vee-bir2006, cos-say2004}, or a Laplace transform of the heat
kernel~\cite{jia-gre-wan2015}.

Rather than taking the Laplace transform of the heat kernel, we take the
approach of using the Laplace transform to convert the time-dependent
parabolic PDE to an parameter-based elliptic PDE. This requires a method
to solve the resulting elliptic PDE, and to invert the Laplace
transform~\cite{LTS2016}. Other groups have also combined BIEs with the
Laplace transform and its inverse including for non-homogeneous heat
conduction in bounded 2D and 3D domains~\cite{sat-pau-gra2002,
sla-sla-zha2003, che-abo-bad1992}.

A complementary approach to solving parabolic and elliptic PDEs, is to
recast their solutions as probability distributions that can be sampled
through Monte Carlo simulation of related diffusion
problems~\cite{MAURO2014, GROSS2022, HWANG20101089, SCHUMM2023,
Lawley2024}. Particle based methods can bypass obstacles inherent to
stencil solution methods, for example resolving singularities at
corners. However, they are slow to converge, relatively low accuracy,
and prone to getting stuck near reflecting
boundaries~\cite{che-lin-her-qua2022}. An additional advantage of the PDE approach developed here is that the density $p(\xx,t)$ is obtained everywhere which yields much more detail on the trajectories of diffusing particles.

In this work, we combine a well-conditioned BIE formulation for the Laplace-transformed PDE, which
involves both Dirichlet and Neumann boundary conditions, and then invert
the Laplace transform with spectral accuracy by integrating along an
appropriate contour in the complex plane. We demonstrate this method for the simulation of \eqref{eqn:diffusion} in the presence of absorbing and reflecting obstacles as shown in the schematics Fig.~\ref{fig:intro}.

The structure of the paper is as follows. In
section~\ref{sec:formulation}, we detail our solution approach based on
a combination of the Laplace transform with boundary integral equations.
In section~\ref{sec:method}, we describe numerical implementation
details, in particular quadrature rules for the integrals and
calculation of surface fluxes. In section~\ref{sec:Numerics}, we provide
a series of examples that both validate the expected convergence rate of
the method and demonstrate its effectiveness on examples with complex
geometries. In particular, our results demonstrate a \emph{shielding}
effect where capture statistics are modulated by geometric features.
Finally in section~\ref{s:conclusions} we conclude by highlighting
avenues for future investigations.




\section{Formulation}
\label{sec:formulation} 
Consider the Laplace transform variable
\begin{align}
  P(\xx,s) = \LL[p](s) = \int_{0}^{\infty} e^{-st} p(\xx,t)\, dt, \quad
  s \in \CC
\end{align}
where a variable with a capital letter denotes the Laplace transform of
the variable with the corresponding lower case letter. Taking the
Laplace transform of~\eqref{eqn:diffusion}, we obtain the elliptic PDE
\bsub
  \label{eqn:diffusionLaplace}
  \begin{alignat}{3}
    (s - \Delta) P(\xx,s) &= \delta(\xx - \xx^*), \qquad 
      &&\xx \in \Omega, \\
    P(\xx,s) &= 0, &&\xx \in \Gamma_D, \\
    \pderiv{P}{\nn}(\xx,s) &= 0, &&\xx \in \Gamma_N.
  \end{alignat} 
\esub
Note that the solution of this PDE depends on $s \in \CC$. For a given
$s$, we write $P$ as
\begin{align}
  \label{eqn:solnrep}
  P(\xx,s) = P^{h}(\xx,s) + G(\xx-\xx^*),
\end{align}
where 
\begin{align}
  \label{eqn:fundsoln}
  G(\xx) = \frac{1}{2\pi} K_0\left( \sqrt{s} |\xx| \right),
\end{align}
is the fundamental solution of the differential operator $s - \Delta$.
Note that the fundamental solution depends on the Laplace transform
variable $s$. To satisfy boundary conditions, we require that
$P^{h}(\xx,s)$ satisfies
\begin{subequations}
  \label{eqn:homoPDE}
  \begin{alignat}{3}
    (s - \Delta) P^{h} &= 0, &&\xx \in \Omega, \\
    P^{h}(\xx,s) &= -G(\xx-\xx^*), &&\xx \in \Gamma_D, \\
    \pderiv{P^{h}}{\nn}(\xx,s) &= -\pderiv{G}{\nn}(\xx-\xx^*), 
      \qquad &&\xx \in \Gamma_N.
  \end{alignat}
\end{subequations}
Having a method to solve equation~\eqref{eqn:diffusionLaplace}, we
compute $p(\xx,t)$ by applying the inverse Laplace transform as
described in Section~\ref{sec:talbot}. We note that since the PDE is
linear, our method can be generalized to the case where the initial
condition is a finite sum of weighted delta functions centered at
different points in $\Omega$. In the following sections we describe
integral equation methods to solve both the homogeneous
PDE~\eqref{eqn:homoPDE} and to calculate the Bromwich integral that
inverts the Laplace transform.

\subsection{Boundary Integral Equation Formulation} 
\label{sec:bies}
We solve the mixed boundary condition elliptic PDE~\eqref{eqn:homoPDE}
using a BIE method. In this manner, we are able to resolve complex
geometries, achieve high-order accuracy, and satisfy the far-field
boundary condition. We account for the mixed boundary conditions by
choosing a layer potential formulation that results in a
well-conditioned system of integral equations. In particular, we express
$P^{h}(\xx,s)$ as a sum of double- and single-layer potentials 
\begin{align}
  \label{eqn:LP}
  P^h(\xx,s) = \DD[\sigma_D](\xx) + \SS[\sigma_N](\xx),
\end{align}
where
\begin{align}
  \DD[\sigma_D](\xx) &= \frac{1}{2\pi} \int_{\Gamma_D} \pderiv{}{\nn_\yy}
    G(\xx - \yy) \sigma_D(\yy)\, ds_\yy, \\
  \SS[\sigma_N](\xx) &= \frac{1}{2\pi} \int_{\Gamma_N} G(\xx - \yy)
    \sigma_N(\yy) \, ds_\yy,
\end{align}
and $G$ is defined in~\eqref{eqn:fundsoln}. We use the unknown density
functions $\sigma_N$ and $\sigma_D$ to satisfy the boundary conditions
by matching the limiting value of $P^h(\xx,s)$ as $\xx$ approaches
$\Gamma_N$ and $\Gamma_D$ with the appropriate boundary condition. In
particular, the density functions satisfy the system of second-kind
Fredholm integral equations
\begin{subequations}
  \label{eqn:BIE}
  \begin{alignat}{3}
    -G(\xx - \xx^*) &= +\frac12 \sigma_D(\xx) + 
      \DD[\sigma_D](\xx) +
      \SS[\sigma_N](\xx), 
      &&\xx \in \Gamma_D, \\
    -\pderiv{}{\nn} G(\xx - \xx^*) &= -\frac12 \sigma_N(\xx) + 
      \pderiv{}{\nn}\DD[\sigma_D](\xx) +
      \pderiv{}{\nn}\SS[\sigma_N](\xx), 
      \qquad &&\xx \in \Gamma_N,
  \end{alignat}
\end{subequations}
where the normal derivatives are taken at the target point $\xx$. The
system of integral equations can be written in matrix form as
\begin{align}
  -\left[
    \begin{array}{c}
      G \\ \partial_\nn G
    \end{array}
  \right] = \left[
    \begin{array}{cc}
      \frac12 & 0 \\ 0 & -\frac12
    \end{array}
  \right] \left[
    \begin{array}{c}
      \sigma_D \\ \sigma_N
    \end{array}
  \right] + \left[ 
    \begin{array}{cc}
      \DD & \SS \\ \partial_\nn \DD & \partial_\nn \SS
    \end{array}
  \right] \left[
    \begin{array}{c}
      \sigma_D \\ \sigma_N
    \end{array}
  \right],
\end{align}
where we have dropped the independent variable. Note that the diagonal
terms of the diagonal blocks are continuous (see Section~\ref{sec:BIE}),
and the off-diagonal blocks are continuous since the target and source
points are on different connected components of $\Gamma$. Therefore,
all the integral operators in~\eqref{eqn:BIE} are compact.

\subsection{Computing the Flux Through $\Gamma_D$}
\label{sec:flux}
To compute the flux $j(t)$, as defined in equation~\eqref{eqn:bdflux},
we first take the normal derivative and integrate the Laplace transform
of $p$. That is, we compute the flux of equation~\eqref{eqn:solnrep},
where $P^h$ is a sum of layer potentials~\eqref{eqn:LP}. In particular,
we compute
\begin{align}
  \label{eqn:LaplaceFlux}
  J(s) = \int_{\Gamma_D} \pderiv{}{\nn} P(\xx,s) \, ds_\xx = 
    \int_{\Gamma_D} \pderiv{}{\nn} \left( P^h(\xx,s) + 
    G(\xx - \xx^*)\right) \, ds_\xx.
\end{align}
We also compute the Laplace transform of the cumulative flux of $j$,
which by properties of the Laplace transform is $C(s) = s^{-1} J(s)$.
Finally, the flux and the cumulative flux are computed in the time
variable by taking an inverse Laplace transform as described in
Section~\ref{sec:talbot}.

\begin{figure}[htbp]
  \centering
  \includegraphics[width=0.45\textwidth]{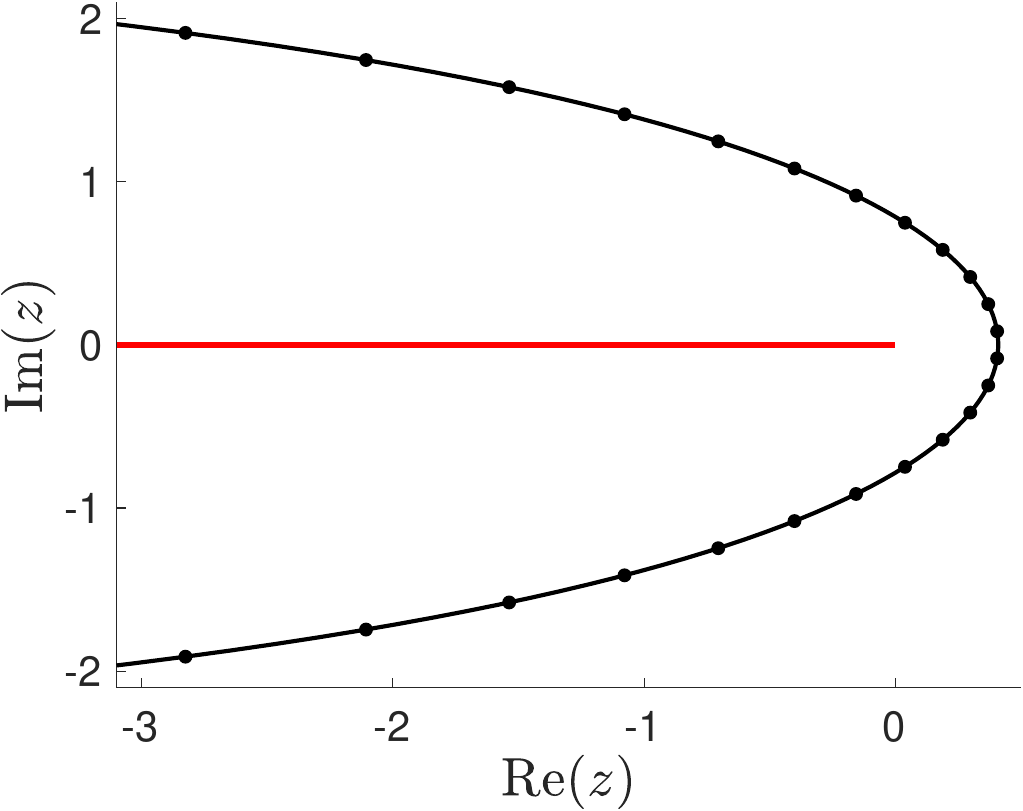}
  \caption{\label{fig:Talbot} The Talbot (solid black)
  contour~\eqref{eqn:talbot}. Integration along the Talbot contour is
  done with the midpoint rule, and the black marks are the resulting
  quadrature points with $M = 24$. The singularities of $P(\xx,s)$ are
  along the negative real axis (solid red).}
\end{figure}

\subsection{Inverting the Laplace Transform}
\label{sec:talbot}
The inverse Laplace transform, also known as the Bromwich integral, is
\begin{align}
  \label{eqn:bromwich}
  \mathcal{L}^{-1}[P](\xx,s) = \frac{1}{2\pi i} 
    \int_{B} e^{st} P(\xx,s)\, ds, \quad \xx \in \Omega,
\end{align}
where
\begin{align}
  B = \{\alpha + i s \:|\: s \in \RR\},
\end{align}
and $\alpha$ is greater than the real part of all the poles of
$P(\xx,s)$. The contour $B$ is called the Bromwich contour. The poles of
$P$ all lie on the negative real axis, so we can choose any value
$\alpha > 0$. Along the Bromwich contour $B$, the integrand $e^{st}
P(\xx,s)$ decays algebraically and oscillates wildly, especially for
large $t$. Therefore, an unpractical number of quadrature points is
required to numerically compute the Bromwich integral. 

Alternatively, we can replace the contour $B$ with any other contour, so
long as all the poles of $P(\xx,s)$ remain on the same side of this new
contour. There are a variety of choices~\cite{tre-wei2014}, and we opt
to use the Talbot contour, 
\begin{align}
  \label{eqn:talbot}
  T = \left\{ \frac{2M}{t}\left(
    -0.6122 + 0.5017\theta\cot(0.6407\theta) + 0.2645i\theta\right)
    \:|\: \theta \in (-\pi,\pi) \right\},
\end{align}
which was first introduced by Dingfelder and
Weideman~\cite{Weideman2015}. The contour depends on both the time $t$
and the number of quadrature points $2M$. We use $2M$ total quadrature
points on the Talbot contour because we are working with real-valued
functions in the time domain, and the integral can be computed along
points on the Talbot contour with $\theta \geq 0$, and then doubling the
result. Although the Talbot contour $T$ does not extend to infinity, the
magnitude of the integrand at $\theta = \pm \pi$ is less than machine
epsilon, so the error introduced by integrating along $T$ rather than
$B$ is numerically zero. In contrast to the Bromwich contour, the
integrand $e^{st} P(s,t)$ is smooth and does not oscillate along the
Talbot contour $T$. The contour is illustrated in
Figure~\ref{fig:Talbot} (black curve), and the poles of $P(\xx,s)$ are
all located along the red line. The black points are quadrature points
that are described in Section~\ref{sec:numericalILT}.

\section{Numerical Methods}\label{sec:method}
Section~\ref{sec:formulation} introduced several integrals that need to
be evaluated in order to solve the heat equation~\eqref{eqn:diffusion}
and the corresponding flux~\eqref{eqn:bdflux}. This section describes
the quadrature techniques we use to solve the second-kind Fredholm
integral equation~\eqref{eqn:BIE} and to evaluate the double-layer
potential~\eqref{eqn:LP} (Section~\ref{sec:BIE}), to compute the flux
through $\Gamma_D$ (Section~\ref{sec:numericalflux}), and to compute the
Talbot integral (Section~\ref{sec:numericalILT}).

\subsection{Solving the BIE}
\label{sec:BIE}
We assume that there are $M_D$ components of $\Gamma_D$ and $M_N$
components of $\Gamma_N$. Each component is parameterized as
$\yy^{k}(\alpha)$, where $k=1,\ldots,M_D + M_N$ and $\alpha \in
[0,2\pi)$. The first $M_D$ components are the Dirichlet components, and
the other $M_N$ components are the Neumann components. Then, the system
of BIEs~\eqref{eqn:BIE} are discretized with a collocation method by
discretizing $[0,2\pi)$ at $N$ equispaced points. We approximate the
integrals with the trapezoid rule, resulting in the linear system
\begin{subequations}
\label{eqn:disBIE}
\begin{align}
  f^k_n &= \frac12 \sigma^k_n + 
    \frac{2\pi}{N}\sum_{m = 1}^{M_D} \sum_{\ell = 1}^{N}
    \pderiv{}{\nn^m_\ell}G\left(\yy^k_n - \yy^m_\ell\right) 
    \sigma^m_\ell + 
    \frac{2\pi}{N}\sum_{m = M_D+1}^{M_D+M_N} \sum_{\ell = 1}^{N}
    G\left(\yy^k_n - \yy^m_\ell\right) \sigma^m_\ell, \nonumber \\
  \label{eqn:disBIE1} 
  &k = 1,\ldots,M_D,\, n = 1,\ldots,N \\
  g^k_n &= -\frac12 \sigma^k_n + 
    \frac{2\pi}{N}\sum_{m = 1}^{M_D} \sum_{\ell = 1}^{N}
    \pderiv{}{\nn^k_n}
    \pderiv{}{\nn^m_\ell}G\left(\yy^k_n - \yy^m_\ell\right) 
    \sigma^m_\ell + 
    \frac{2\pi}{N}\sum_{m = M_D+1}^{M_D+M_N} \sum_{\ell = 1}^{N}
    \pderiv{}{\nn^k_n}
    G\left(\yy^k_n - \yy^m_\ell\right) \sigma^m_\ell, \nonumber \\
  \label{eqn:disBIE2} 
  &k = M_D+1,\ldots,M_D+M_N,\, n = 1,\ldots,N,
\end{align}
\end{subequations}
where 
\begin{align}
  \sigma_{\ell}^m = \sigma\left(\yy^{m}\left(\alpha_\ell\right)\right)
    \left|\pderiv{}{\alpha}\yy^{m}\left(\alpha_\ell\right)\right|,
\end{align}
and the left-hand sides are
\begin{align}
  f^k_n = -G\left(\yy^k_n - \yy^*\right), \text{ and } 
  g^k_n = -\pderiv{}{\nn^k_n}G\left(\yy^k_n - \yy^*\right).
\end{align}
Our notation uses a superscript to denote the body being considered, and
the subscript denotes the discretization point on that particular body.
Since the superscript determines if the density function is defined on a
Neumann or Dirichlet body, we do not need to specify whether the
discretized version of corresponds to $\sigma_D$ or $\sigma_N$.

Because the fundamental solution involves the modified Bessel function
of the second kind, the second summand in equation~\eqref{eqn:disBIE1}
has a logarithmic singularity when $\ell = n$ and $m = k$. However, this
never occurs since $k$ and $m$ never coincide. For the same reason, the
first summand in equation~\eqref{eqn:disBIE2} has no singularities.
However, the first summand of equation~\eqref{eqn:disBIE1}, and the
second summand of equation~\eqref{eqn:disBIE2} cannot be naively
evaluated when $\ell = n$ and $m = k$, but the singularities of the
kernels are removable. Therefore, the terms in the summation when $\ell
= n$ and $m = k$ are replaced with their limiting values
\begin{align}
  \lim_{\substack{\yy \rightarrow \xx \\ \xx \in \Gamma}} 
    \pderiv{}{\nn_\yy} G(\xx - \yy) = -\frac{1}{2}\kappa(\xx),
    \qquad
  \lim_{\substack{\yy \rightarrow \xx \\ \xx \in \Gamma}} 
    \pderiv{}{\nn_\xx} G(\xx - \yy) = +\frac{1}{2}\kappa(\xx),
\end{align}
where $\kappa(\xx)$ is the curvature at the surface point
$\xx\in\Gamma$. 

A numerical issue with BIE solvers is near-singular integration which
occurs when two different components of $\Gamma$ are close to one
another. Independent of the value of $N$, the error of the trapezoid
rule will grow without bound. There are a suite of quadrature methods
that address nearly-singular integrals, but these are not the focus of
this work. Instead, we use geometries whose components are sufficiently
separated that the problems can achieve sufficient accuracy for values
of $N$ that can be implemented on a standard laptop. The smoothness of
the kernels, particularly the regularity of the kernels of the
double-layer potential and normal derivative of the single-layer
potential, determines the convergence with respect to $N$. Since these
kernels are both normal derivatives of the modified Bessel function, the
kernel is $C^1$, and therefore the trapezoid rule achieves third-order
accuracy~\cite{kro-qua2011}. This spatial convergence behavior is
demonstrated in Section~\ref{sec:BIEconvergence}. We note that
higher-order accuracy is possible using a variety of quadrature
methods~\cite{hel2013, alp1999, kap-rok1997}.

Since we are solving a second-kind Fredholm integral equation, and all
the kernels are continuous, the discretization~\eqref{eqn:disBIE} can be
solved with a mesh-independent number of GMRES
iterations~\cite{cam-ips-kel-mey-xue1996}. Therefore, the computational
cost of solving~\eqref{eqn:disBIE} is proportional to the cost of
computing the necessary summations at all target points. We use a direct
calculation to compute these sums, therefore requiring $\mathcal{O}((M_D
+ M_N)^2 N^2)$ operations. In future work, we will use a direct solver
or a fast summation method to reduce the computational cost.

\subsection{Computing the Flux}
\label{sec:numericalflux}
The Laplace transform of the flux into the absorbing bodies is given by
equation~\eqref{eqn:LaplaceFlux}. We must compute the normal derivative,
followed by a boundary integral around $\Gamma_D$, of two terms: the
Green's function centered at $\xx^{*}$, and the layer potential
representation of $P^h(\xx,s)$. The former of these two terms is
computed exactly by calculating the normal derivative of
equation~\eqref{eqn:fundsoln}, and since $\xx^{*} \notin \Gamma_D$, the
integral is computed with spectral accuracy with the trapezoid rule.
However, the latter term requires quadrature to compute
\begin{align}
  \int_{\Gamma_D}& \pderiv{}{\nn} P^h(\xx,s) \, ds_\xx 
  =\int_{\Gamma_D} \pderiv{}{\nn} 
    \left(\DD[\sigma_D](\xx) + \SS[\sigma_N](\xx) \right) ds_\xx
    \nonumber \\
  &=\int_{\Gamma_D} \pderiv{}{\nn}\left(
    \int_{\Gamma_D} \pderiv{}{\nn_\yy}G(\xx - \yy) \sigma_{D}(\yy) \, ds_\yy + 
    \int_{\Gamma_N} G(\xx - \yy) \sigma_{N}(\yy) \, ds_\yy 
    \right) \, ds_\xx 
    \nonumber \\
  &=\int_{\Gamma_D} \int_{\Gamma_D} 
    \frac{\partial^2}{\partial \nn \partial \nn_\yy}G(\xx - \yy)
    \sigma_{D}(\yy) \, ds_\yy ds_\xx + 
    \int_{\Gamma_D} \int_{\Gamma_N} \pderiv{}{\nn}G(\xx - \yy) 
    \sigma_{N}(\yy) \, ds_\yy ds_\xx.
  \label{eqn:fluxLP}
\end{align}
Since $\xx \in \Gamma_D$ and $\yy \in \Gamma_N$, and these two contours
do not intersect, we accurately compute the second integral in equation~\eqref{eqn:fluxLP} with the trapezoid rule. That is
\begin{align}
  \int_{\Gamma_D} \int_{\Gamma_N} \pderiv{}{\nn}G(\xx - \yy) 
    \sigma_{N}(\yy) \, ds_\yy ds_\xx 
  \approx \left(\frac{2\pi}{N}\right)^2
  \sum_{k=1}^{M_D}\sum_{n=1}^{N}
  \sum_{m=M_D+1}^{M_D+M_N}\sum_{\ell=1}^{N} 
  \pderiv{}{\nn_{n}^{k}}
  G\left(\yy_{n}^{k} - \yy_{\ell}^{m}\right) \sigma_{\ell}^m.
\end{align}
The first integral in equation~\eqref{eqn:fluxLP} requires more care
since $\xx \in \Gamma_D$ and $\yy \in \Gamma_D$, and the singularity of
the kernel scales as $|\xx - \yy|^{-2}$ as $\xx \rightarrow \yy$.
Therefore, the trapezoid rule will not converge to the principal value
integral. Instead, the integral is manipulated with a careful use of
singularity subtraction~\cite{nit2020}
\begin{align}
  \int_{\Gamma_D} \int_{\Gamma_D} 
    &\frac{\partial^2}{\partial \nn \partial \nn_\yy}G(\xx - \yy)
    \sigma_{D}(\yy) \, ds_\yy ds_\xx 
  = \frac{1}{2\pi} \int_{\Gamma_D} \int_{\Gamma_D} 
    \frac{\partial^2}{\partial \nn \partial \nn_\yy}
    K_0\left(\sqrt{s}|\xx - \yy| \right)
    \sigma_{D}(\yy) \, ds_\yy ds_\xx \nonumber \\
  &= \frac{1}{2\pi} \int_{\Gamma_D} \int_{\Gamma_D}
    \frac{\partial^2}{\partial \nn \partial \nn_\yy}
    \left(K_0\left(\sqrt{s}|\xx - \yy| \right) + \log|\xx - \yy|\right)
    \sigma_{D}(\yy) \, ds_\yy ds_\xx \nonumber \\
  &-\frac{1}{2\pi} \int_{\Gamma_D} \int_{\Gamma_D}
    \frac{\partial^2}{\partial \nn \partial \nn_\yy} \log|\xx - \yy| 
    \sigma_{D}(\yy) \, ds_\yy ds_\xx.
  \label{eqn:fluxSingularitySubtraction}
\end{align}
Equation~\eqref{eqn:fluxSingularitySubtraction} simply introduces the
addition and subtraction of a term. The particular choice guarantees
that the two integrals can be individually computed accurately or
analytically. Based on the asymptotics of the modified Bessel function,
the singularity of the first integrand scales with $|\xx
- \yy|^{-1}$ as $\xx \rightarrow \yy$, and the integrand is computed
using odd-even integration~\cite{sid-isr1988}. The error of this
quadrature method is first-order accurate, but, as mentioned earlier,
higher-order accuracy is possible by using other more expensive
quadrature rules. The second integral is identically zero since
\begin{align}
  \frac{1}{2\pi} \int_{\Gamma_D} \int_{\Gamma_D}
    \frac{\partial^2}{\partial \nn \partial \nn_\yy} &\log|\xx - \yy| 
    \sigma_{D}(\yy) \, ds_\yy ds_\xx = 
  \frac{1}{2\pi} \int_{\Gamma_D} \pderiv{}{\nn} \int_{\Gamma_D}
    \pderiv{}{\nn_\yy} \log|\xx - \yy| 
    \sigma_{D}(\yy) \, ds_\yy ds_\xx \nonumber \\
  &= \frac{1}{2\pi} \int_{\Omega} \Delta \left(\int_{\Gamma_D}
    \pderiv{}{\nn_\yy} \log|\xx - \yy| 
    \sigma_{D}(\yy) \, ds_\yy\right) d\xx = 0,
\end{align}
where the inner integral is zero because the integrand is the
double-layer potential of Laplace's equation.

\begin{figure}[htbp]
  \centering
  \includegraphics[width=0.4\textwidth]{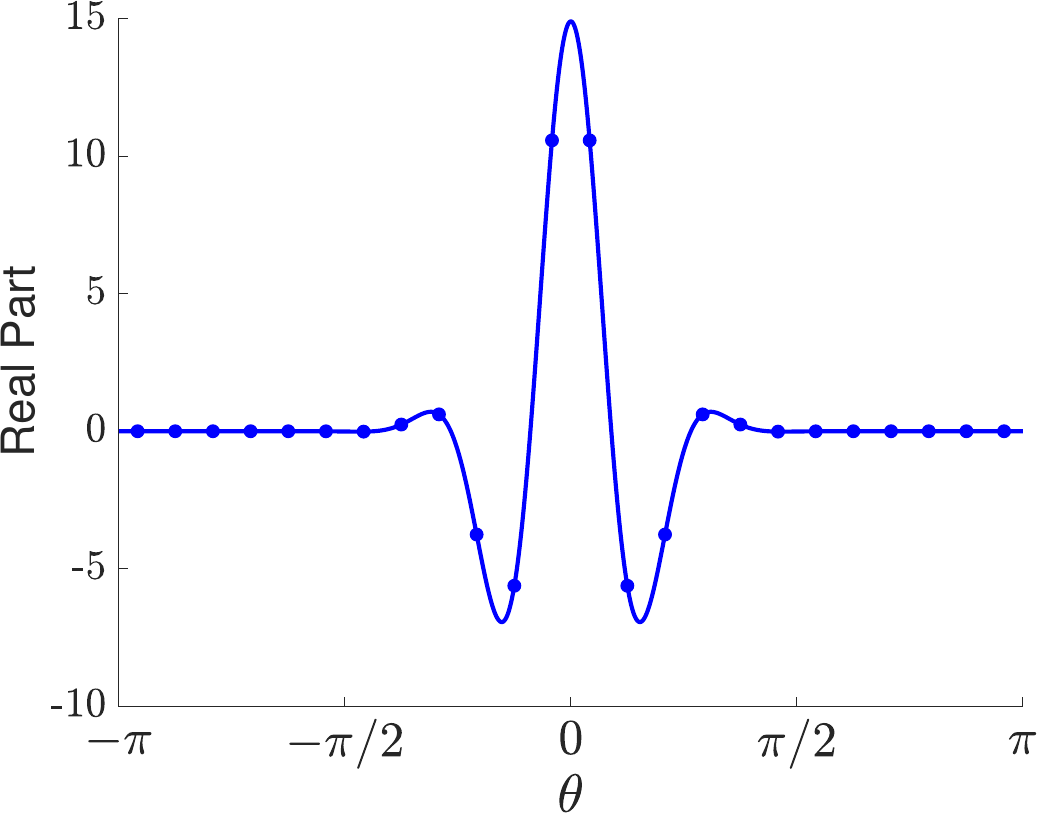}
  \qquad
  \includegraphics[width=0.4\textwidth]{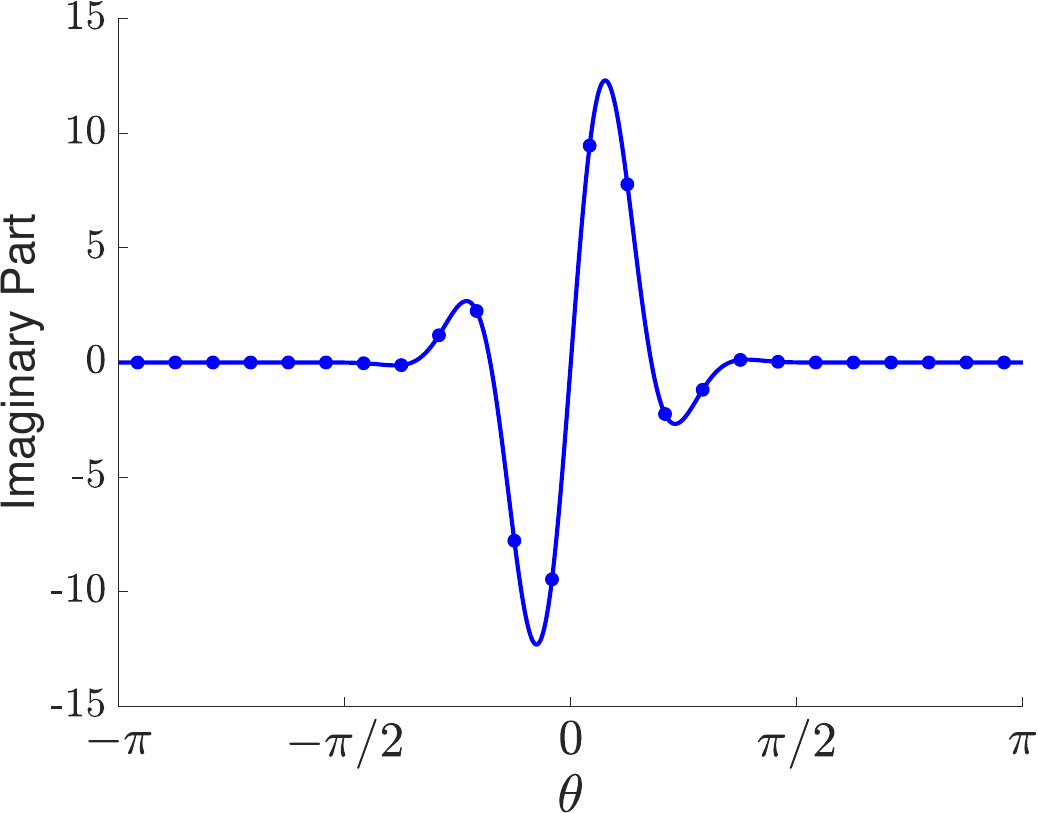}
  \caption{\label{fig:talbotIntegrand} The real and imaginary parts of a
  typical integrand required to compute a Talbot integral. The
  truncation of the integrand and the choice of the Talbot contour
  guarantee that the quadrature error is
  $\mathcal{O}\left(10^{-1.2M}\right)$, until $M$ is sufficiently large
  that machine precision is reached.}
\end{figure}

\subsection{Computing the Inverse Laplace Transform} 
\label{sec:numericalILT}
Figure~\ref{fig:Talbot} illustrates the Talbot contour we use, which, to
the best of our knowledge, achieves the highest-order accuracy known to
date. Spectral accuracy is achieved by applying the midpoint rule, and
the midpoint quadrature points are the black points in
Figure~\ref{fig:Talbot}. The integrand for an example from
Section~\ref{sec:talbot_convergence} is shown in
Figure~\ref{fig:talbotIntegrand}. The marks indicate the quadrature
points used for the trapezoid rule with $M = 12$ and $t=10$. The
magnitude of the integrand at the final quadrature points is on the
order of $10^{-14}$, so truncation errors can be safely ignored.

The lack of oscillatory behavior and the decay rate of the integrand
guarantee that the midpoint rule achieves spectral accuracy with an
error that scales as $10^{-1.2M}$. In
Section~\ref{sec:talbot_convergence}, we demonstrate this temporal order
of convergence by applying the quadrature to the Laplace transform of
the one-dimensional heat equation.

As noted by Dingfeld and Weideman~\cite{Weideman2015}, floating-point
error dominates the midpoint rule if $M$ is too large. They control this
error by adjusting the Talbot contour when necessary. We take this same
approach for the one-dimensional heat equation in
Section~\ref{sec:talbot_convergence}, and observe that this adjustment
should be made when $M > 12$. Therefore, to compute the inverse Laplace
transform for all other numerical examples, we apply the midpoint rule
to equation~\eqref{eqn:talbot} with $M = 12$ quadrature points.

\section{Numerical Results}\label{sec:Numerics}
In this section, we demonstrate the accuracy and utility of the
formulation and numerical methods described in
Sections~\ref{sec:formulation} and~\ref{sec:method}. We first
demonstrate that the expected orders of accuracy are achieved in time
(Section~\ref{sec:talbot_convergence}) and in space
(Section~\ref{sec:BIEconvergence}). Then, we consider several complex
geometries to demonstrate the effect of screening the absorbing body
with a collection of absorbing and reflecting bodies
(Sections~\ref{sec:spiral}--\ref{sec:maze}).

\subsection{Talbot Integral Convergence}
\label{sec:talbot_convergence}
\begin{figure}[htbp]
\begin{center}
  \includegraphics[width=0.45\textwidth]{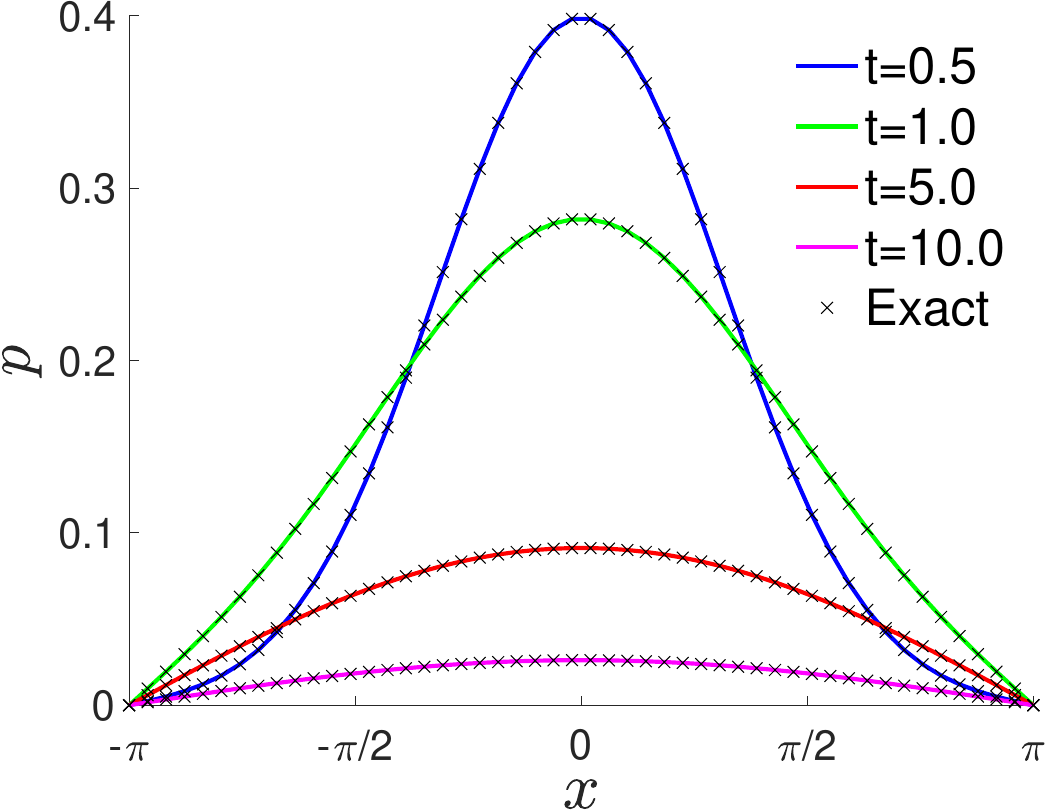} \qquad
  \includegraphics[width=0.45\textwidth]{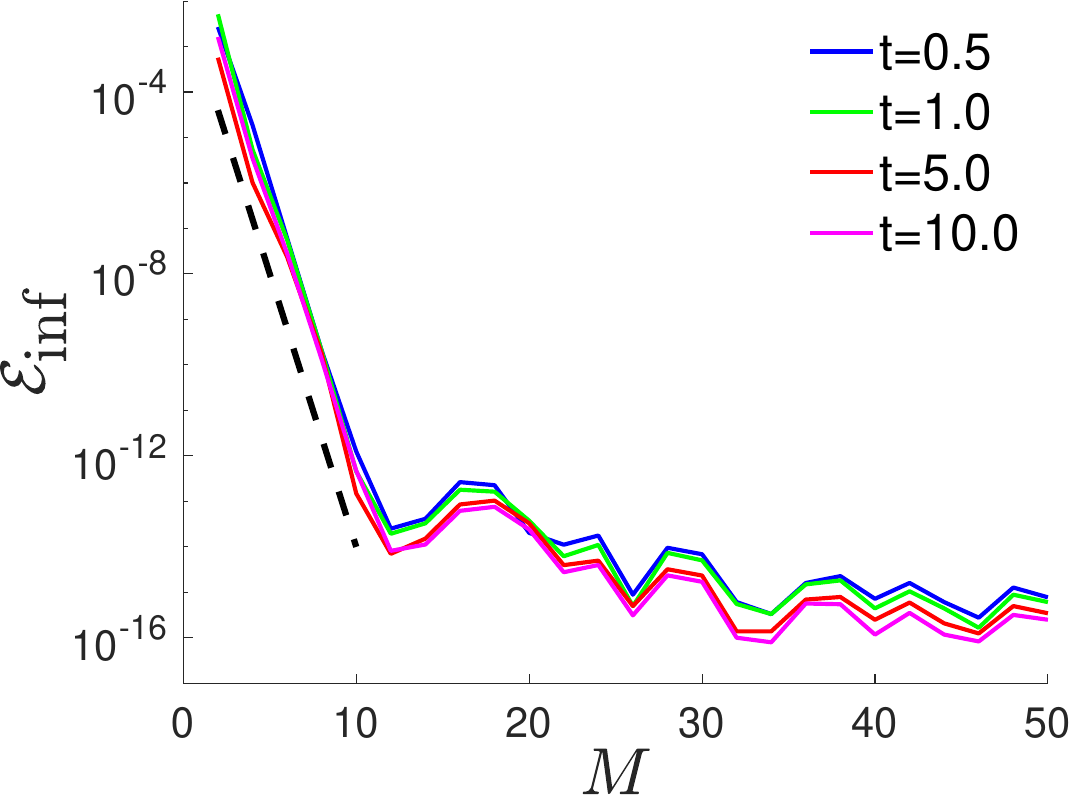}
\end{center}
\caption{\label{fig:1DSolution} A verification of the convergence of
Talbot integration on the exactly solvable 1D heat
equation~\eqref{eqn:1DPeriodicDiffusion}. The left plot shows the exact
and approximate solutions, and the right plot shows the error at four
different times as a function of the number of Talbot quadrature points.
For $M \leq 12$, the error converges as $10^{-1.2 M}$ (black dashed
curve). If the same Talbot contour were used for larger values of $M$,
round-off error would begin to dominate the quadrature error. Therefore,
for these larger values of $M$, the Talbot contour is slightly
modified~\cite{Weideman2015}.}
\end{figure}

We first demonstrate that the method achieves the expected rate of
convergence in time by considering the exactly solvable one-dimensional
heat equation in the interval $x \in [-\pi,\pi]$. Specifically, we solve
\begin{subequations}
  \label{eqn:1DPeriodicDiffusion}
  \begin{alignat}{3}
    \pderiv{p}{t} &= \frac{\partial^2 p}{\partial x^2}, 
    &&x \in (-\pi,\pi), \quad t>0, \\
    p(-\pi,t) &= p(\pi,t) = 0, \qquad && t>0, \\
    p(x,0) &= \delta(x), &&x \in (-\pi,\pi).
  \end{alignat}
\end{subequations}
The Laplace transform $P(x,s) = \LL[p](s) = \int_{0}^{\infty} e^{-st}
p(x,t)\, dt$ solves the modified Helmholtz problem
\bsub
  \label{eqn:1DPeriodicDiffusionLaplace}
  \begin{align}
    sP - \frac{\partial^2 P}{\partial x^2} & = \delta(x), 
      \qquad x \in (-\pi,\pi), \quad s \in \CC, \\
    P(-\pi,s) &= P(\pi,s) = 0,
  \end{align}
\esub
whose exact solution is
\begin{align}
  P(x,s) = \frac{\sinh [\sqrt{s}(\pi-|x|)]}
    {2\sqrt{s}\cosh \sqrt{s} \pi}.
\end{align}
Therefore the solution $p(x,t)$ can be expressed with the inverse
Laplace transform
\begin{align}\label{eq:LaplaceFull}
  p(x,t) = \frac{1}{2\pi i} \int_{B} e^{st} P(x,s) \, ds = 
  \frac{1}{2\pi i} \int_{B} e^{st} \frac{\sinh 
    [\sqrt{s}(\pi-|x|)]}{2\sqrt{s}\cosh \sqrt{s} \pi}\, ds,
\end{align}
where $B$ is the Bromwich contour. Applying the transformation $s = z^2$ yields
\begin{align}
\label{eqn:1DSolution}
  p(x,t) = \frac{1}{2\pi i} \int_{B'} f(z)\, dz, \qquad 
  f(z) =  e^{z^2t} \frac{\sinh[ z(\pi-|x|)]}{\cosh z \pi},
\end{align}
where $B'$ is the transformed Bromwich contour, and this integral can be
evaluated with standard residue calculus. In particular, the integrand
of~\eqref{eqn:1DSolution} has simple poles $z=
i\left(n-\tfrac{1}{2}\right)$ for $n = 1,2,3,\ldots$. To apply the
residue theorem, we identify that
\begin{align}
  \mbox{Res}(f(z), i(n-\tfrac{1}{2})) = \frac{1}{\pi} 
  e^{-\left(n-\tfrac12\right)^2t} \cos \left[(n-\tfrac12)x\right].
\end{align}
Hence we arrive at the solution of~\eqref{eqn:1DPeriodicDiffusion}
\begin{align}\label{eq:SeriesFull}
  p_{\textrm{ex}}(x,t) = 2\pi i \frac{1}{2\pi i}\sum_{n=1}^{\infty}   
  \mbox{Res}\left(f(z), i\left(n-\tfrac{1}{2}\right)\right)  
    = \frac{1}{\pi} \sum_{n=1}^{\infty} 
  e^{-\left(n-\tfrac12\right)^2t} \cos
  \left[\left(n-\tfrac12\right)x\right],
\end{align}
which coincides with the solution derived by separation of variables. We
now evaluate the error of the numerical inverse procedure applied
to~\eqref{eq:LaplaceFull} by comparing with the series
solution~\eqref{eq:SeriesFull}. In particular, the inverse Laplace
transform is computed by applying the $M$-point midpoint rule along the
Talbot contour $T$ defined in equation~\eqref{eqn:talbot}. As the number
of quadrature points increases, round-off error associated with
evaluation of $e^{st}$ becomes sizable. We implement a round-off control
algorithm~\cite[Sec.~3]{Weideman2015} to prevent this issue from
introducing large errors. For a selection of times,
Figure~\ref{fig:1DSolution} shows the exact solution, the numerical
solution, and the error
\begin{align}
  \label{eq:errorInf}
  \mathcal{E}_{\textrm{inf}} =\max_{x\in(-\pi,\pi)} 
    |p(x,t) - p_{\textrm{ex}}(x,t)|,
\end{align}
for different number $M$ of Talbot quadrature points. The expected
exponential order of convergence is achieved until the error reaches
machine precision.

\subsection{Boundary Integral Equation Convergence}
\label{sec:BIEconvergence}
Having validated the convergence of the Talbot integration, we examine
the convergence behavior of the BIE solver by considering a
two-dimensional problem with a closed form solution. We let $\Omega =
\RR^{2} \backslash B(0,1)$ where $B(0,1)$ is the unit ball centered at
the origin, and impose a homogeneous Dirichlet boundary condition on
$\partial\Omega$. The initial condition is the delta function centered
at $\xx^{*} = (x_0,0)$, for $x_0 > 1$.

\begin{figure}[htbp]
\centering
\includegraphics[width=0.4\textwidth]{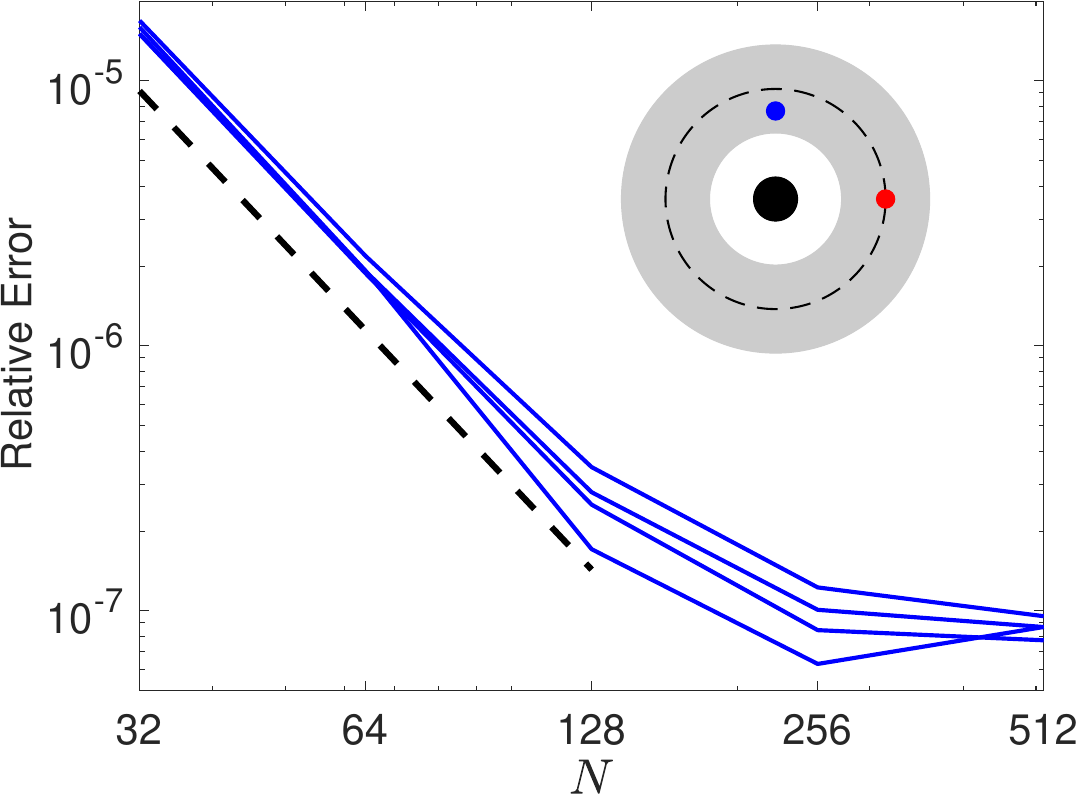}
\includegraphics[width=0.4\textwidth]{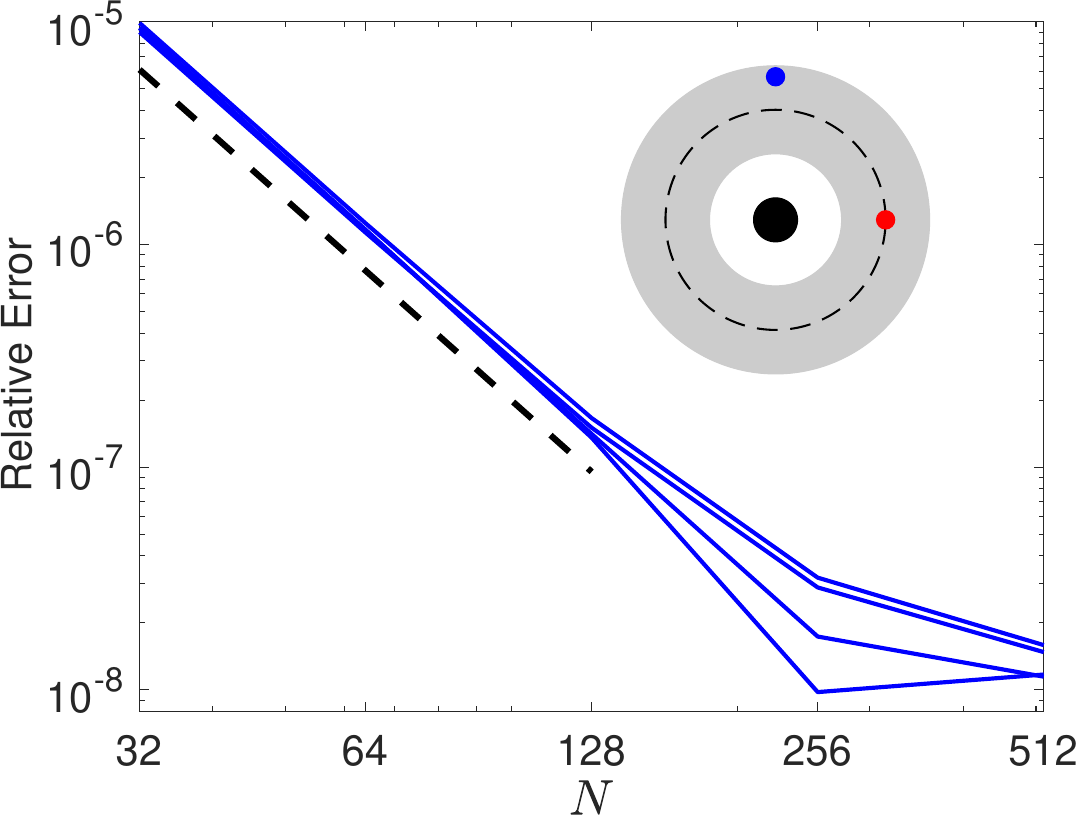}
\includegraphics[width=0.4\textwidth]{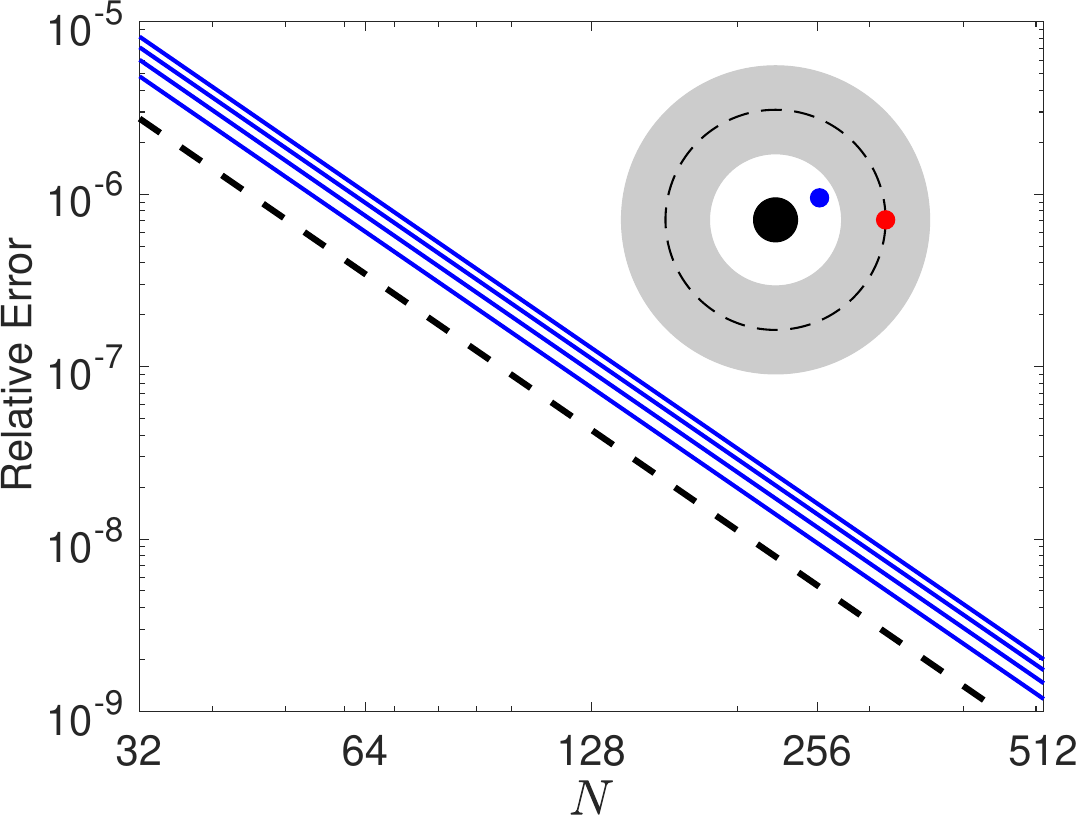}
\includegraphics[width=0.4\textwidth]{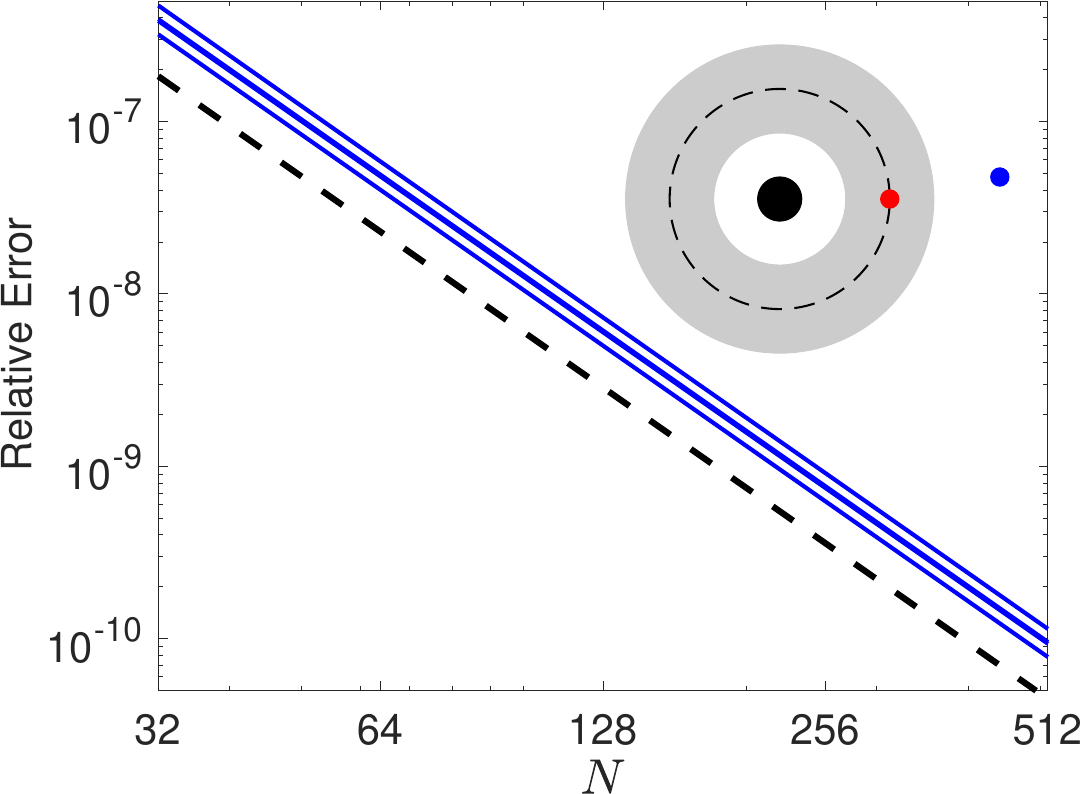} 
\caption{A spatial convergence study of the BIE solver at four target
points for increasing number of collocation points $N$. The black circle
in the middle of the inset is the absorbing body. Each blue curve
represents a particular $s$-value on the Talbot contour at time $t=10$.
The dashed black lines indicate third-order convergence. If the modulus
of the target (blue point) is too close (gray region) to the modulus of
the initial particle location (red point), then the exact solution has
round-off error which results in the loss of convergence for large $N$.
\label{fig:2DConvergence}}
\end{figure}

The governing equation for $P = P(r,\theta;s)$ in polar coordinates $\xx
= r(\cos\theta, \sin \theta)$ with $R = |\xx^{*}|$ is
\bsub\label{eq:heatDisc}
\begin{alignat}{3}
\label{eq:heatDisc_a}   \frac{\partial^2 P}{\partial r^2} +\frac{1}{r}\frac{\partial P}{\partial r} +\frac{1}{r^2}\frac{\partial^2 P}{\partial \theta^2} -sP &=-\frac{1}{r}\delta(r-R)\delta(\theta),
    \qquad &&r>1, \quad \theta \in(-\pi,\pi), \quad s\in\CC,\\[5pt]
\label{eq:heatDisc_b}   P(r,\theta+2\pi) &= P(r,\theta), &&r>1,\\[5pt]
\label{eq:heatDisc_c}   P(1,\theta) &= 0, &&\theta\in(-\pi,\pi).
\end{alignat}
\esub
The exact solution of $P$ is the Fourier cosine
series~\cite[equation~(A5)]{ste-zmu-hen-rij-ors-lin2020}
\bsub\label{eqn:C_sol}
\begin{align}\label{eqn:C_sol_a}
P(r,\theta;s)= 
\begin{cases}
\ds\sum_{n=0}^{\infty} A_n\left[ I_n(\sqrt{s} r) -
  \ds\frac{I_n(\sqrt{s})}{K_n(\sqrt{s})}K_n(\sqrt{s}r) \right]\cos
  n\theta, & r \in (1,R],\\[10pt]
\ds\sum_{n=0}^{\infty} A_n\left[ \frac{I_n(\sqrt{s}R)}{K_n(\sqrt{s}R)} -
  \frac{I_n(\sqrt{s})}{K_n(\sqrt{s})} \right]K_n(\sqrt{s}r) \cos n\theta, & r\in(R,\infty),
\end{cases}
\end{align}
where the constants $A_n$ are
\begin{align}\label{eqn:C_sol_b}
A_n = \begin{cases} 
  \frac{1}{2\pi} K_0(\sqrt{s}R), & n = 0,\\[5pt] 
  \frac{1}{\pi} K_n(\sqrt{s}R), & n \geq 1. 
\end{cases}
\end{align}
\esub
The series solution~\eqref{eqn:C_sol} suffers from overflow and
underflow error when $|\xx|$ and $|\xx^{*}|$ are too close, and this is
demonstrated in Figure~\ref{fig:2DConvergence}. In particular, we
consider four different locations $\xx$ where the error of the numerical
solution is computed for increasing values of $N$. If $|\xx|$ is
sufficiently smaller or larger than $|\xx^*|$, then the expected
third-order convergence is observed for different values of $s$.
However, when $|\xx|$ and $|\xx^*|$ are too close (gray region),
underflow results in an incorrect exact solution, and ultimately results
in the error plateauing above machine precision as $N$ increases.

In addition, we can calculate a series solution for the Laplace
transform of the flux,
\begin{align}
J(s) = \int_{\Gamma_D} P_r|_{_{r=1}} ds =   
    \frac{K_0(\sqrt{s}R)}{K_{0}(\sqrt{s})}.
\end{align}
Since the Laplace transform of the cumulative flux satisfies $sC(s) =
J(s)$, we obtain the series solution for $C(s)$,
\begin{align}
  \label{eqn:exactOneBodyFlux}
  C(s) = \frac{J(s)}{s} = \frac{K_0(\sqrt{s}R)}{sK_{0}(\sqrt{s})}.
\end{align}

\subsection{Navigating a Complex Geometry}
\label{sec:spiral}
\begin{figure}[htbp]
\centering
\includegraphics[height=0.38\textwidth]{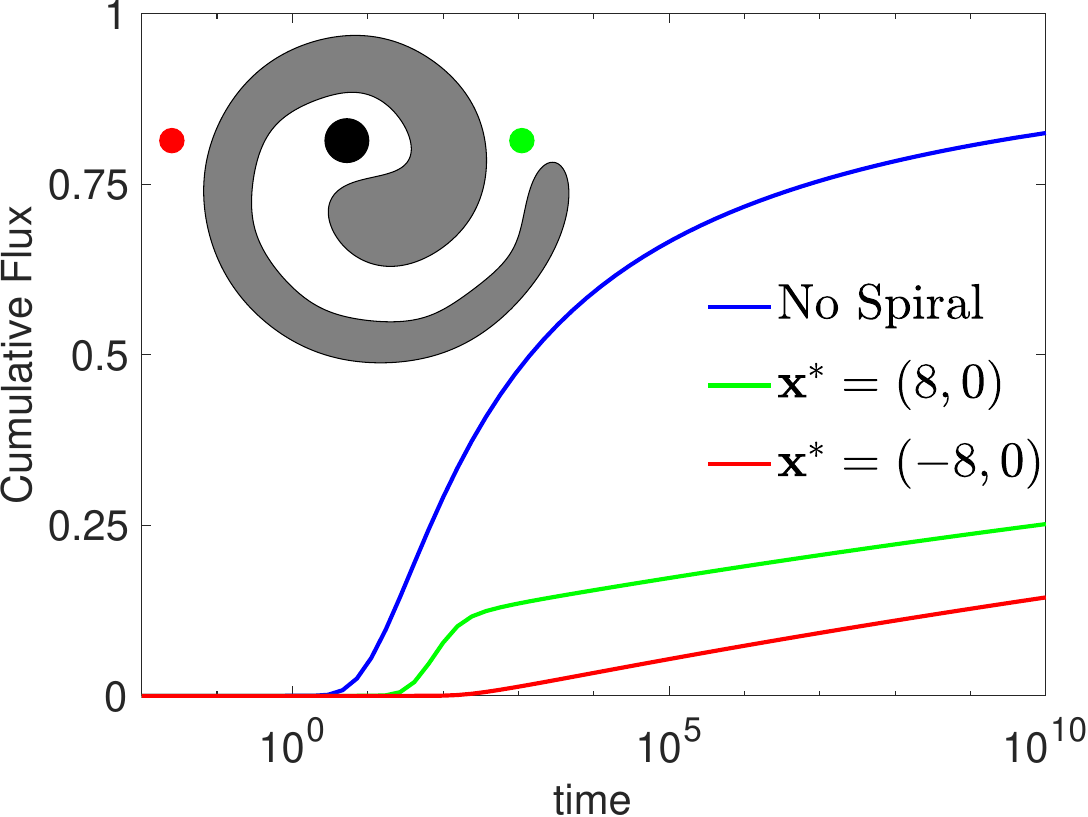}
\caption{\label{fig:spiralCDF} Three cumulative flux into the absorbing
trap (black circle) in the middle of a reflecting spiral (gray body).
Two experiments are initiated by particles at $\xx^* = (-8,0)$ (red) and
$\xx^* = (8,0)$ (green). The blue curve displays the cumulative flux in
the absence of the spiral and hence the capture rate in the absence of
geometric effects.}
\end{figure}

We consider a particle navigating to a single absorbing body in the
middle of a single {\em spiral-shaped} reflecting body
(Figure~\ref{fig:spiralCDF}). The absorbing body is a unit circle
centered at the origin, and we solve the governing equations with two
different sources: $\xx^* = (8,0)$ is near the opening of the spiral
(green mark), and $\xx^* = (-8,0)$ is at the back of the spiral (red
mark). Since the initial particle location of the latter example is
farther from the opening of the spiral, we expect later arrival times
into the absorbing body. As expected, the cumulative flux remains nearly
zero for a very long time, then slowly starts to increase around
$t=10^2$. With the spiral reflecting body, the probability of capture by
$t = 10^{10}$ from $\xx^* = (8,0)$ is only $25\%$ while from $\xx^* =
(-8,0)$, the equivalent probability is only $15\%$. As a baseline, we
include the cumulative flux when the reflecting spiral is absent (blue
curve).

\begin{figure}[htp]
\centering
\includegraphics[width=0.303\textwidth,trim=0pt 0pt 75pt
0pt,clip=true]{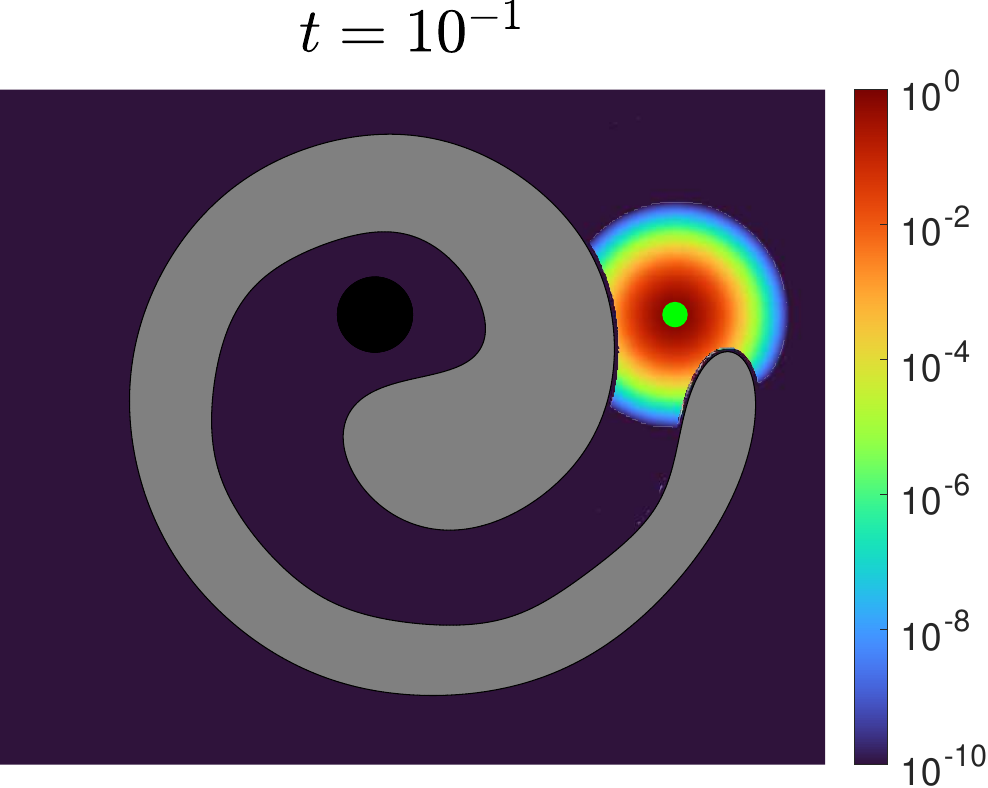}
\includegraphics[width=0.303\textwidth,trim=0pt 0pt 75pt
0pt,clip=true]{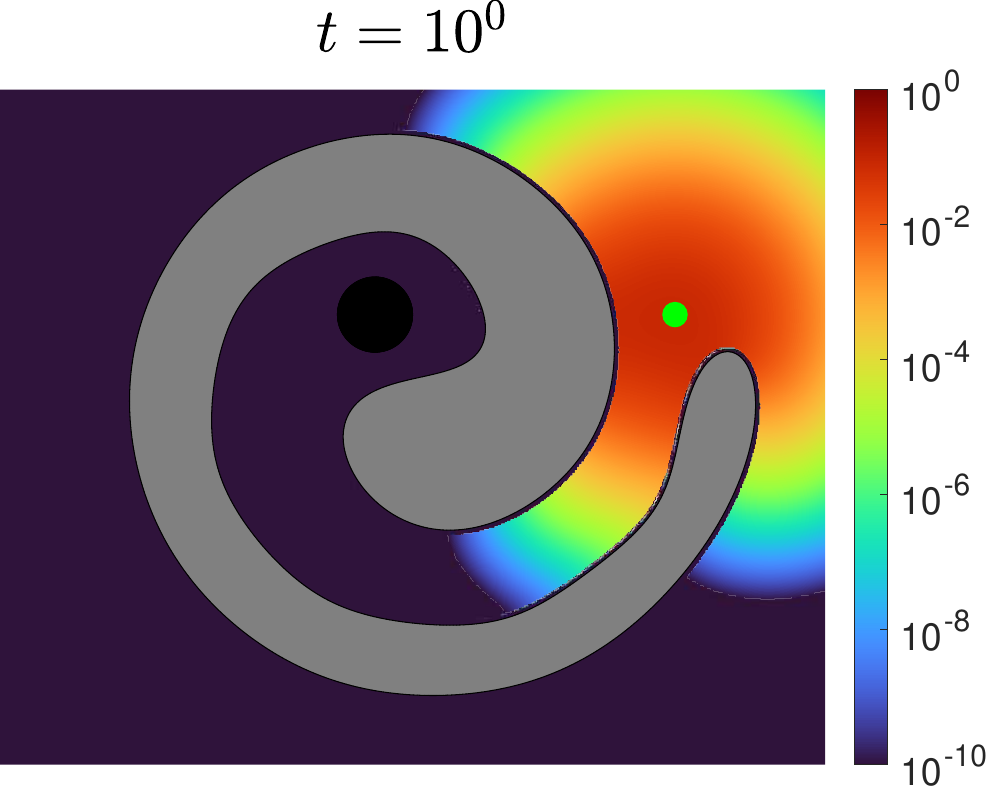}
\includegraphics[width=0.36\textwidth]{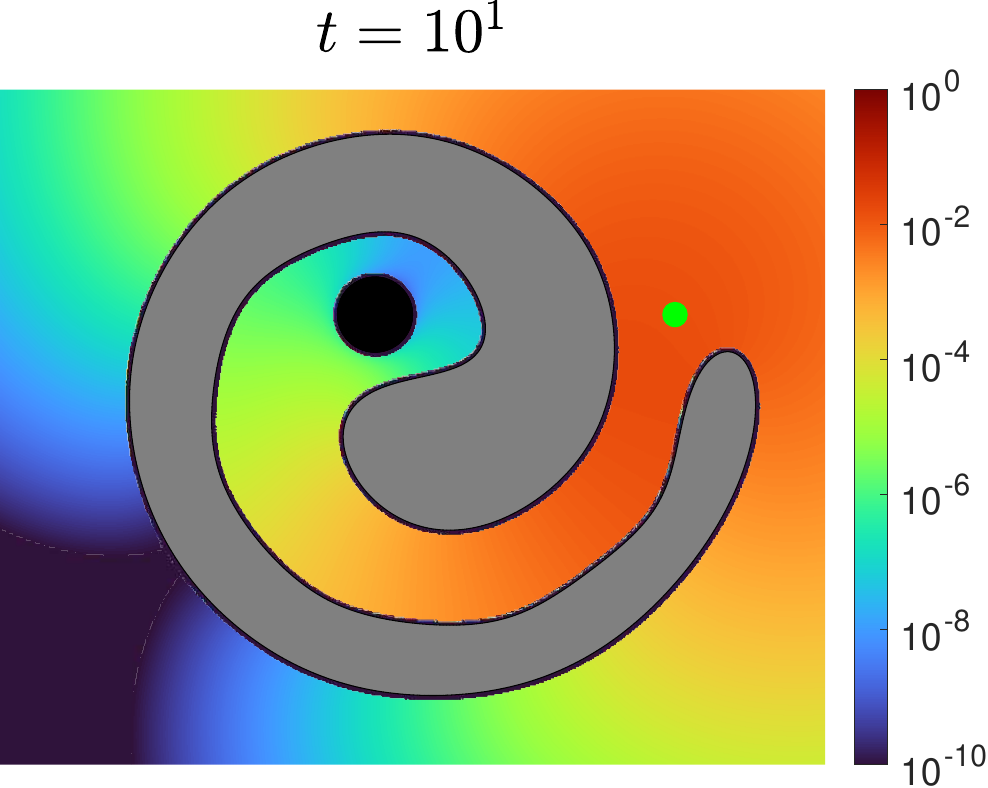}
\\[10pt]

\includegraphics[width=0.303\textwidth,trim=0pt 0pt 75pt
0pt,clip=true]{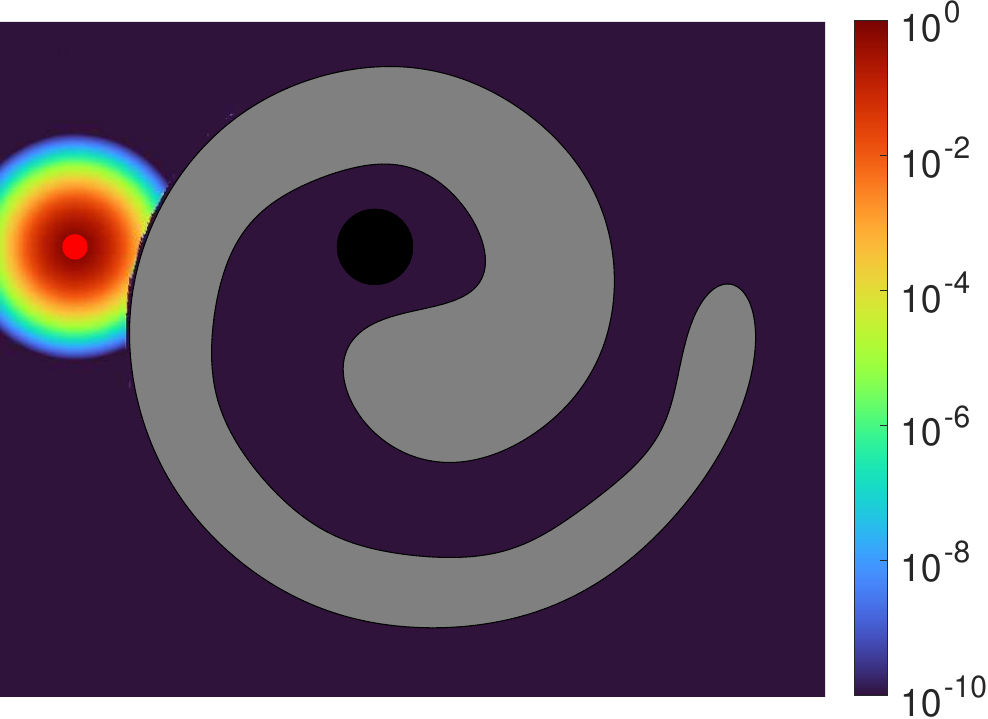}
\includegraphics[width=0.303\textwidth,trim=0pt 0pt 75pt
0pt,clip=true]{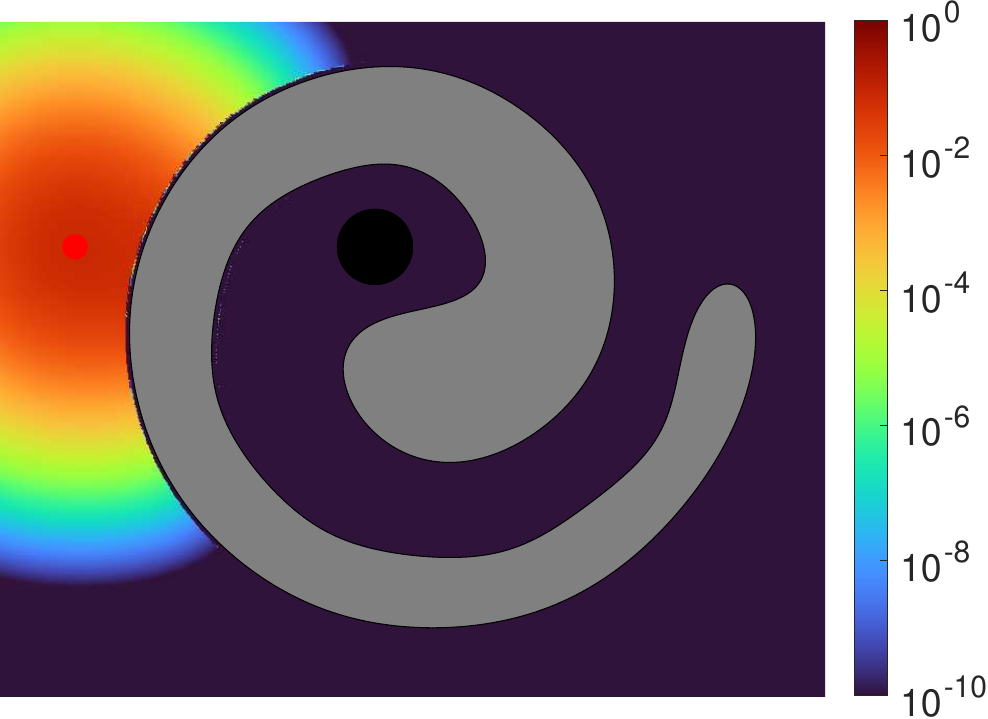}
\includegraphics[width=0.36\textwidth]{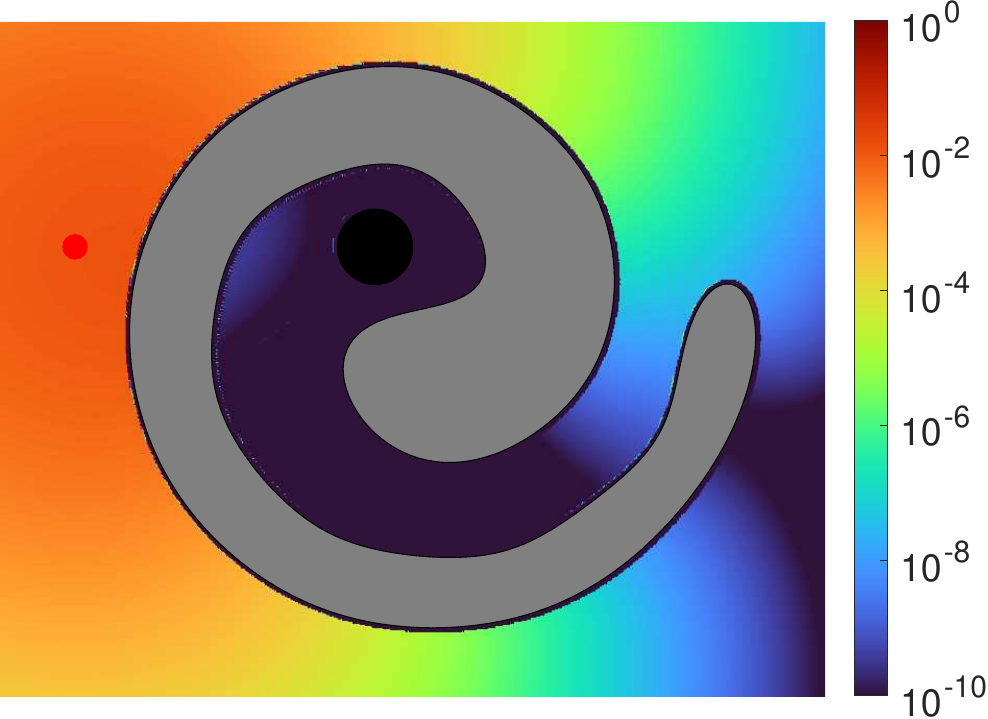}
\caption{\label{fig:spiralHeat} A heat map of solutions $p(\xx,t)$
of~\eqref{eqn:diffusion} at three times $t =\{10^{-1},10^0,10^1\}$. The
gray spiral is reflecting, the black disk is absorbing, and the particle
is initialized at the green dot (top row) and red dot (bottom row).}
\end{figure}

A heat map of $p(\xx,t)$ is shown in Figure~\ref{fig:spiralHeat} for
both initial locations, and at three values of $t$. Consistent with the
cumulative flux, we observe that a particle initiated near the opening
of the spiral is much more likely to be found near the absorbing body
when compared to a particle initiated at the location opposite the
spiral opening. Hence, this demonstrates a simple geometric scenario
where the capture rate is highly modulated by geometric features.

\subsection{Shielding Effects in the Faraday Cage}
\label{sec:faraday}
\begin{figure}[htp]
\centering
\includegraphics[width=0.5\textwidth]{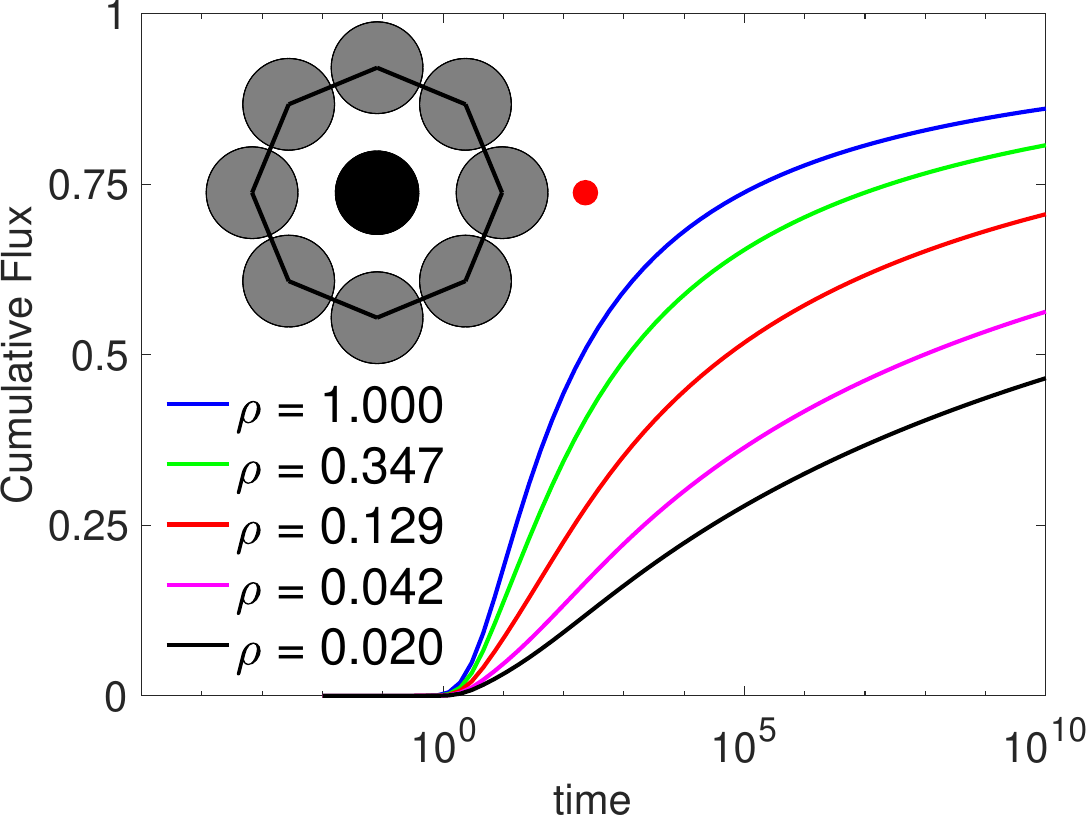}
\caption{\label{fig:faradayCDF} The cumulative flux into an absorbing
trap (black circle) in the middle of eight reflecting bodies (gray
bodies). The particles are initiated at the red dot. The confining ratio
$\rho$ is the ratio of the perimeter of the octagon (black curve) that
is outside the reflecting bodies relative to the perimeter of the entire
octagon. The confining ratio of the inset is $\rho = 0.129$.}
\end{figure}

Inspired by the Faraday cage, which blocks electromagnetic fields from
an enclosure~\cite{Faraday2015}, we investigate the effect of shielding
a Brownian particle by placing eight reflecting bodies around an
absorbing body. The absorbing body is a unit circle centered at the
origin, and the particle source is initially positioned at
$\xx^*=(5,0)$. The centers of the reflecting bodies are evenly
distributed at the eight roots of unity on a circle of radius three. We
increase the shielding effect by increasing the radii of these
reflecting bodies. This setup, with reflecting bodies having a radius of
$1.1$, is illustrated in the inset of Figure~\ref{fig:faradayCDF}.
Given a particular radius, we define a confining number ($\rho$) by
considering the octagon passing through the center of the reflecting
bodies and define $\rho$ to be the ratio of the octagon's perimeter that
is outside of the reflecting bodies to the total perimeter of the
octagon.

The cumulative flux into the absorbing body is shown in
Figure~\ref{fig:faradayCDF}. Each line corresponds to a different
confining ratio $\rho$. The $\rho = 1$ curve corresponds to the scenario
with no reflecting bodies, and the exact solution in Laplace space is given by equation~\eqref{eqn:exactOneBodyFlux}.
 The case $\rho = 0$ corresponds
to the vanishing gap limit in which the absorber is not accessible from
the initial location. As expected, increasing the radii of the
reflecting bodies (lowering $\rho$) decreases the cumulative flux at any
particular time. In the most shielded example ($\rho = 0.020$), not even
50\% of the cumulative flux of particles are absorbed at time $t =
10^{10}$. These distributions have so much weight in their tails that
Monte Carlo methods are not appropriate for these examples.

\begin{figure}[htp]
\centering
\includegraphics[height=0.27\textwidth, 
  trim=0pt 0pt 75pt 0pt,clip=true]{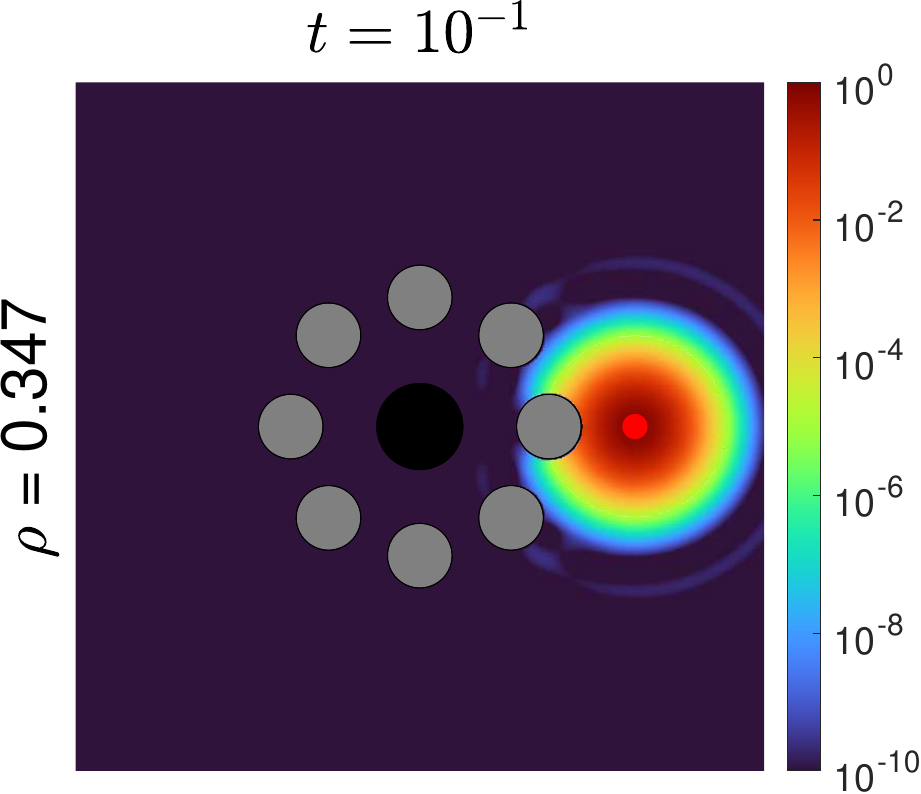}
\includegraphics[height=0.27\textwidth, 
  trim=0pt 0pt 75pt 0pt,clip=true]{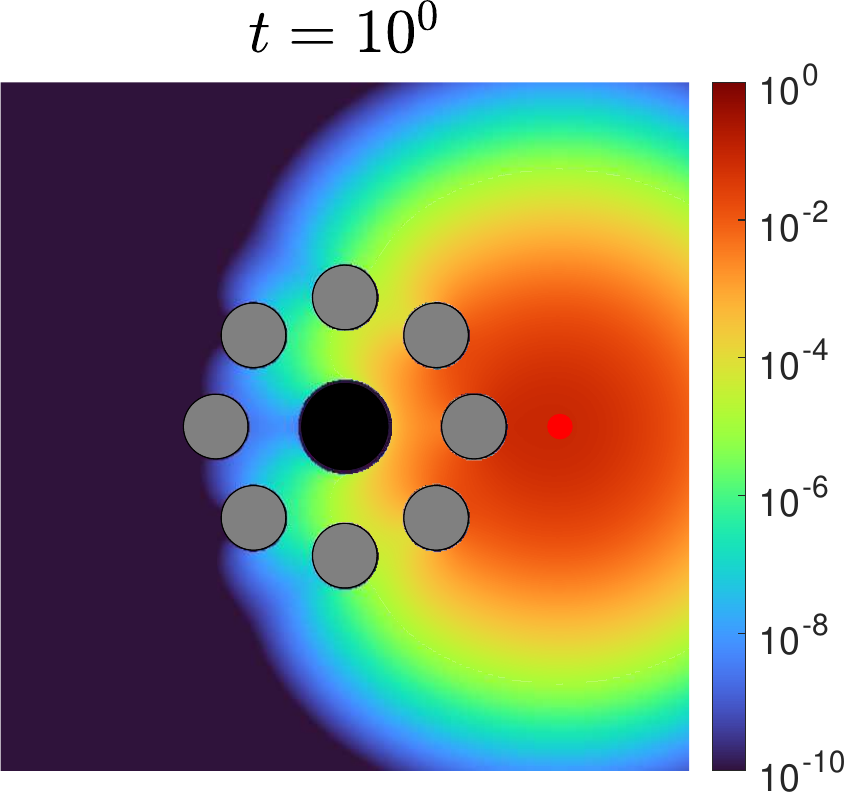}
\includegraphics[height=0.27\textwidth]{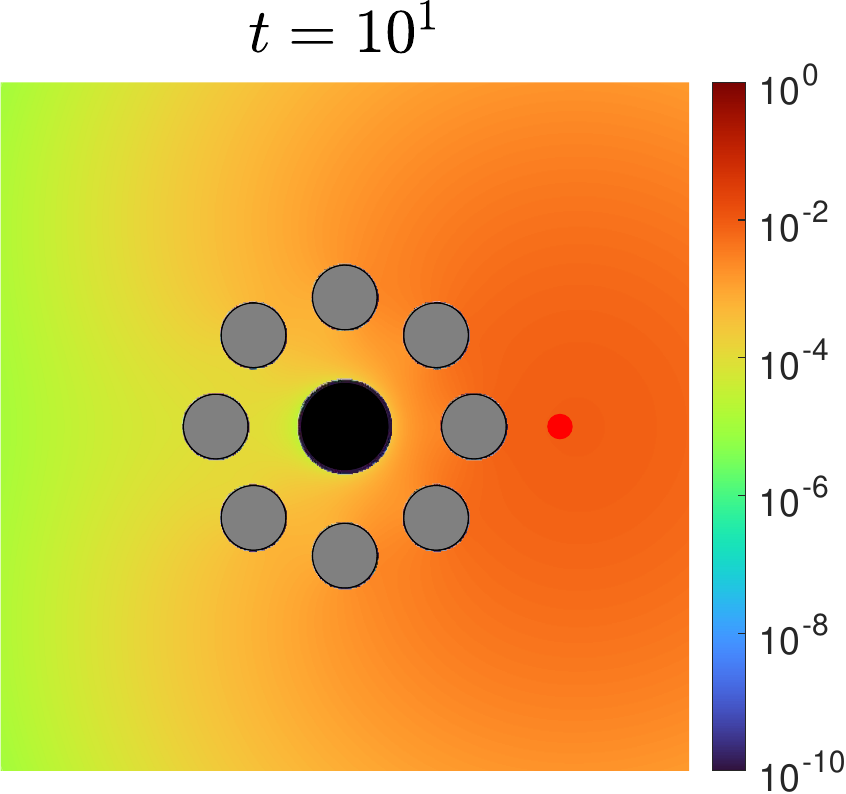}
\\[2pt]
\includegraphics[height=0.2485\textwidth, 
  trim=0pt 0pt 75pt 0pt,clip=true]{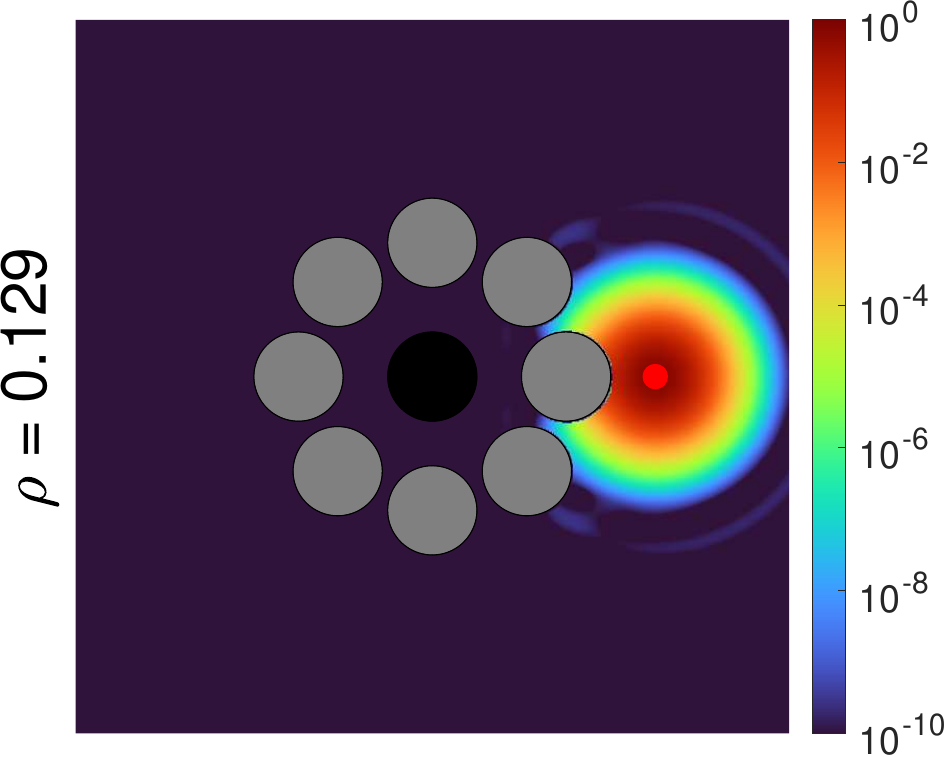}
\includegraphics[height=0.2485\textwidth, 
  trim=0pt 0pt 75pt 0pt,clip=true]{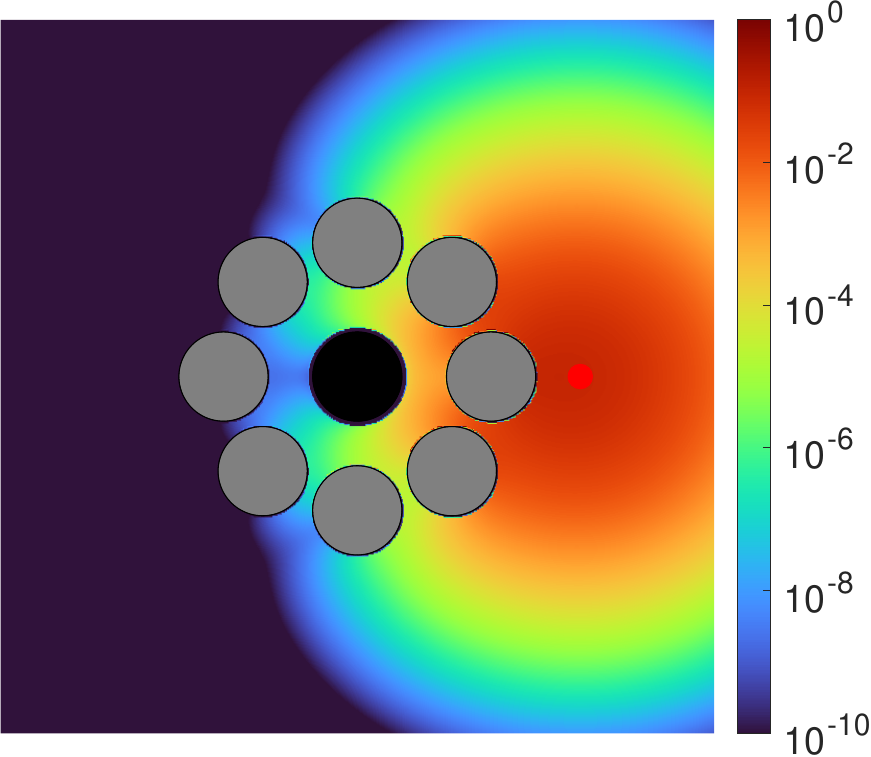}
\includegraphics[height=0.2485\textwidth]{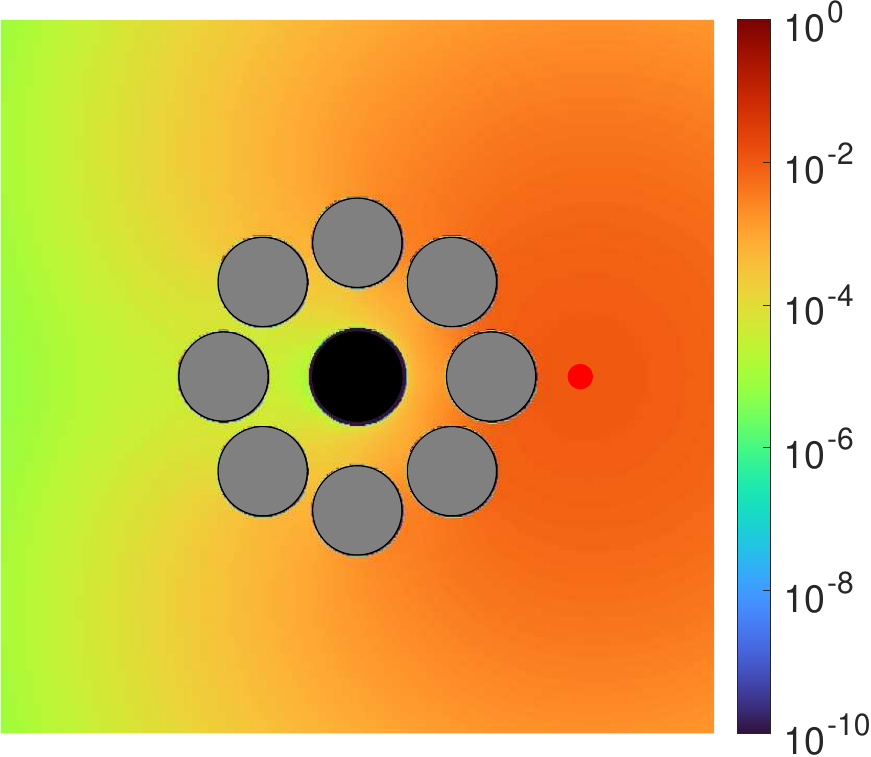}
\\[2pt]
\includegraphics[height=0.2485\textwidth, 
  trim=0pt 0pt 75pt 0pt,clip=true]{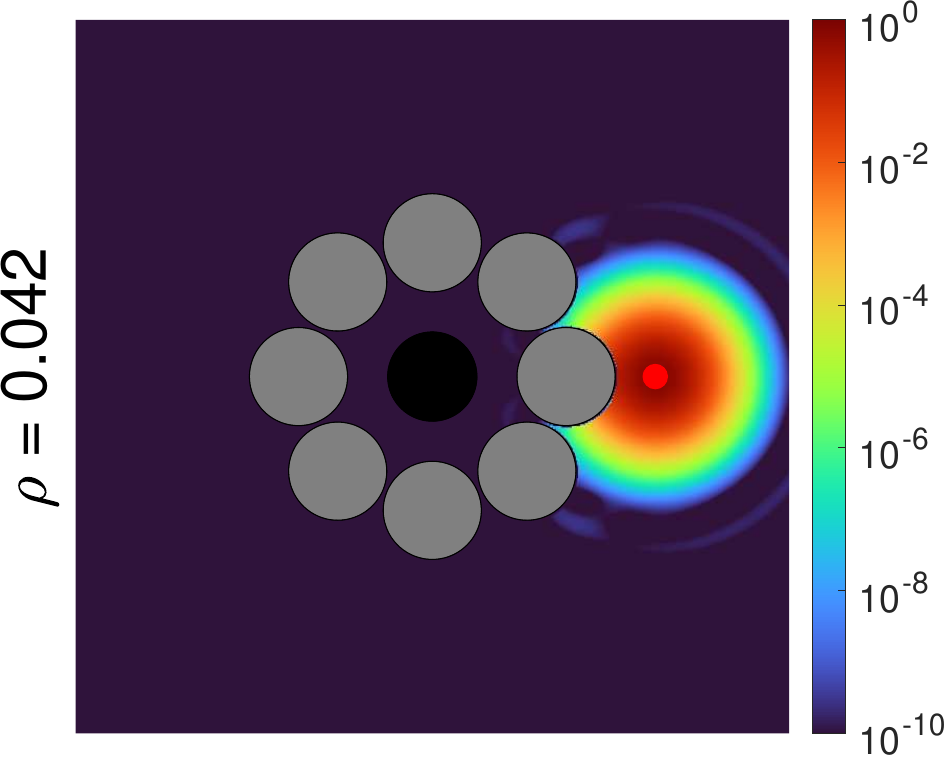}
\includegraphics[height=0.2485\textwidth, 
  trim=0pt 0pt 75pt 0pt,clip=true]{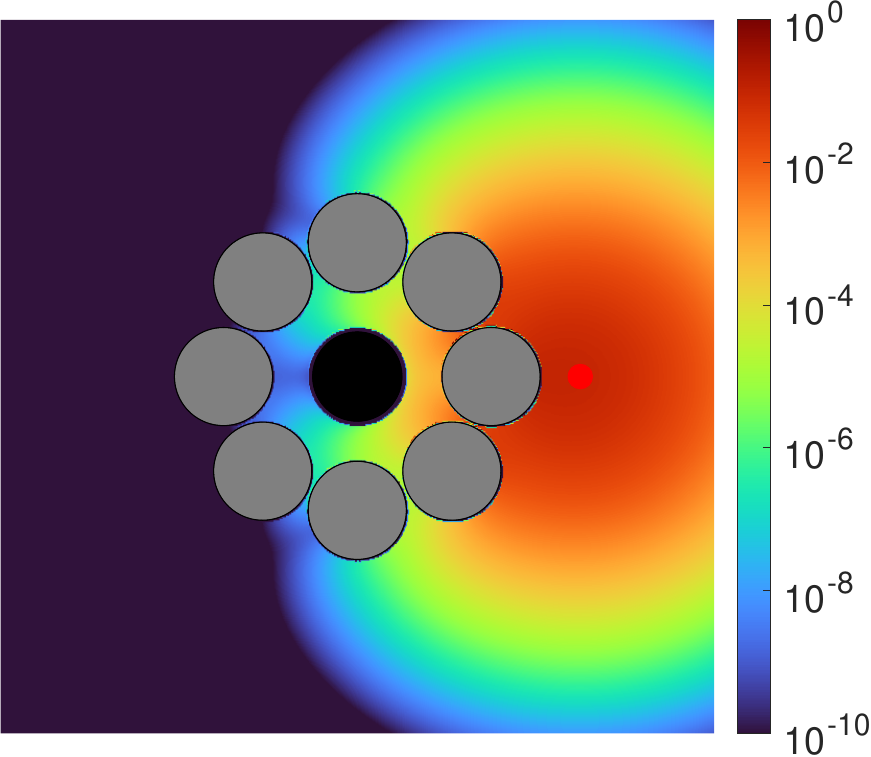}
\includegraphics[height=0.2485\textwidth]{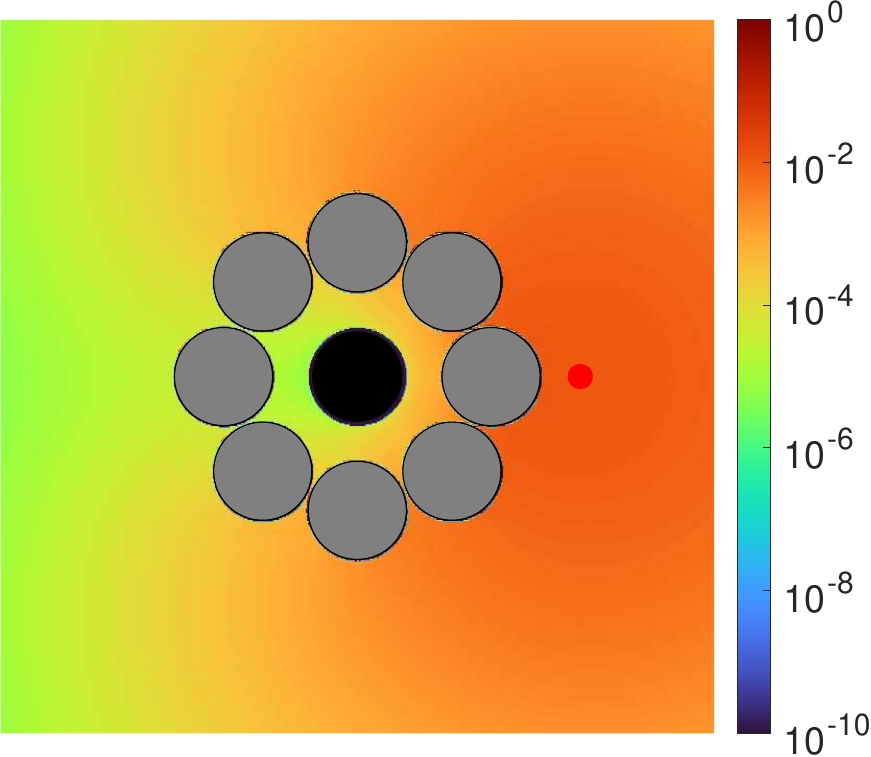}
\\[2pt]
\includegraphics[height=0.2485\textwidth, 
  trim=0pt 0pt 75pt 0pt,clip=true]{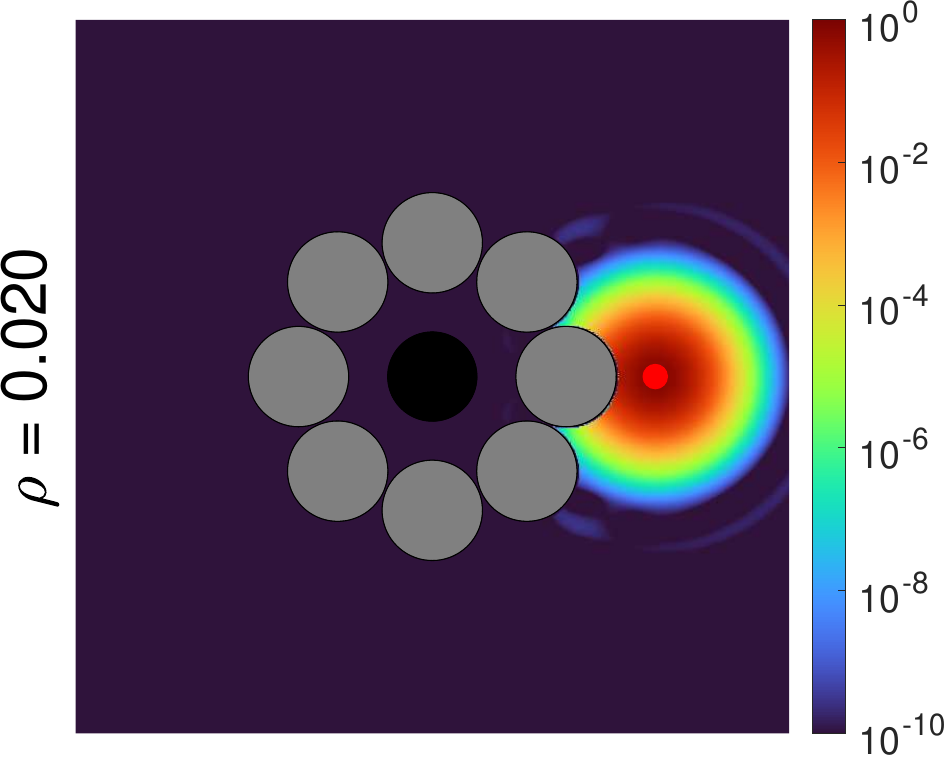}
\includegraphics[height=0.2485\textwidth, 
  trim=0pt 0pt 75pt 0pt,clip=true]{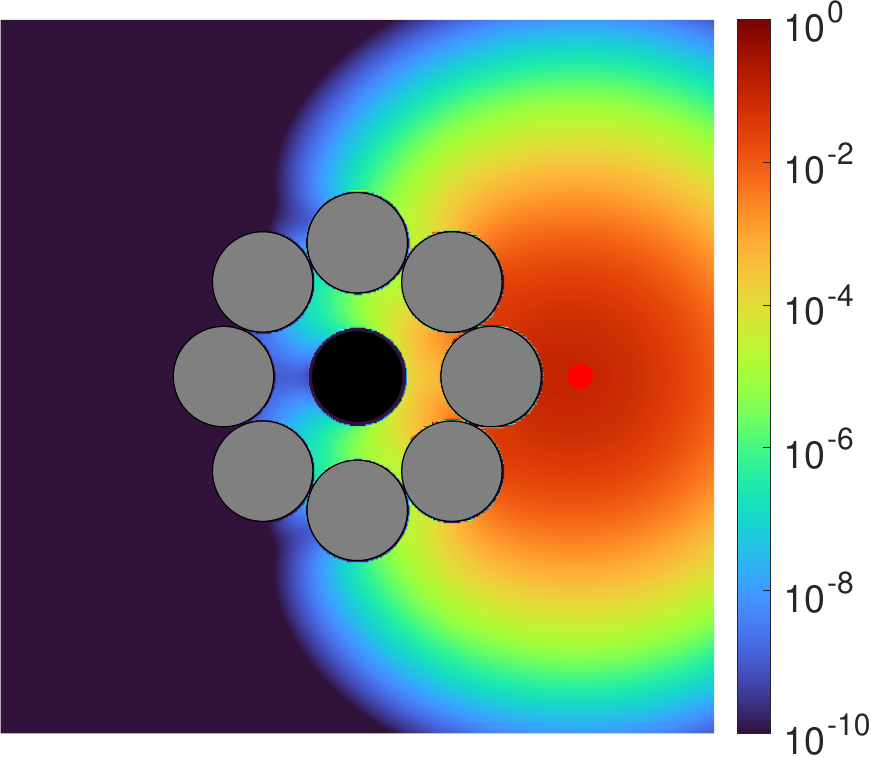}
\includegraphics[height=0.2485\textwidth]{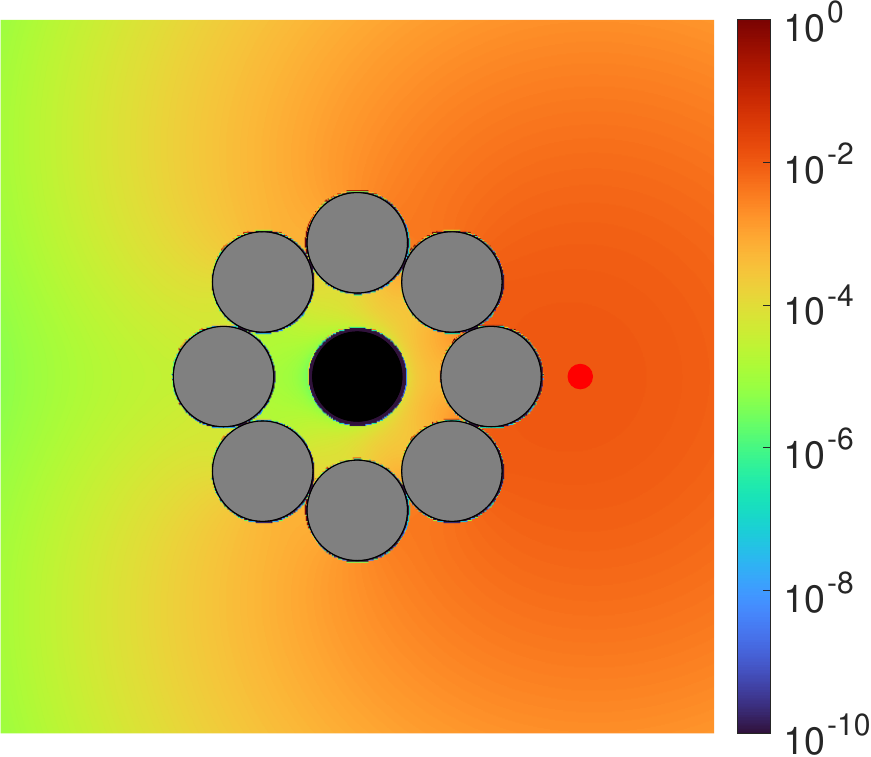}
\caption{\label{fig:faradayHeat} A heat map of solutions $p(\xx,t)$ to \eqref{eqn:diffusion} at short times $t = \{10^{-1},10^0,10^1\}$. The gray circles are reflecting, the black circle is absorbing,
and the particle is initialized at the red dot. From top to bottom, the geometries are defined by their confining ratios $\rho = \{0.347,0.129,0.042,0.020 \}$.}
\end{figure}

\begin{figure}[htp]
\centering
\includegraphics[height=0.23\textwidth,trim=0pt 0pt 75pt
0pt, clip=true]{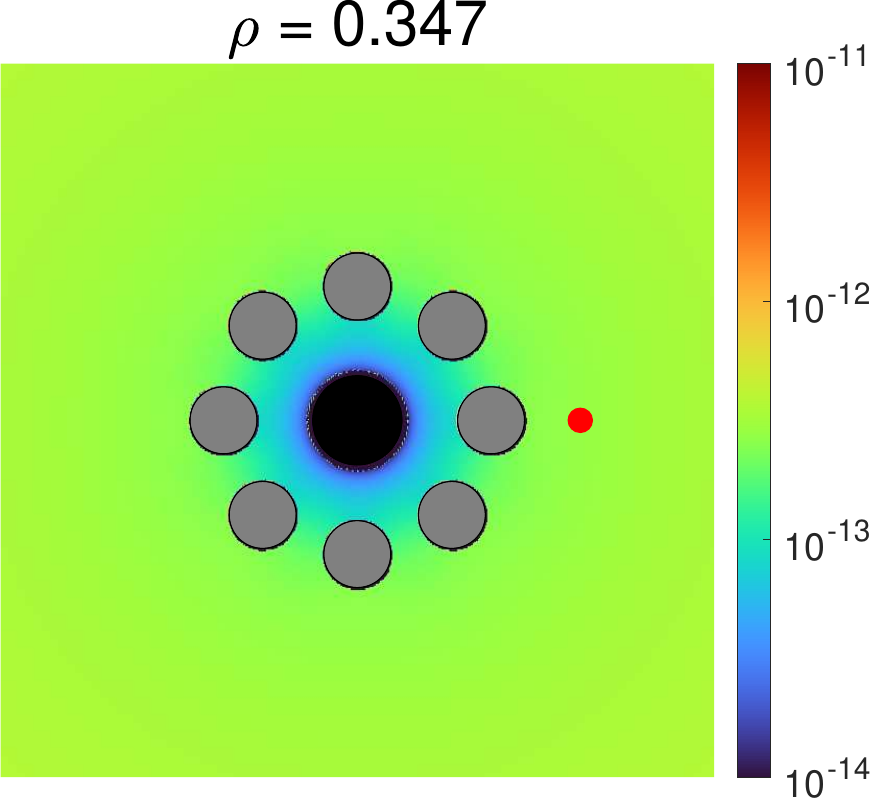}
\includegraphics[height=0.23\textwidth,trim=0pt 0pt 75pt
0pt, clip=true]{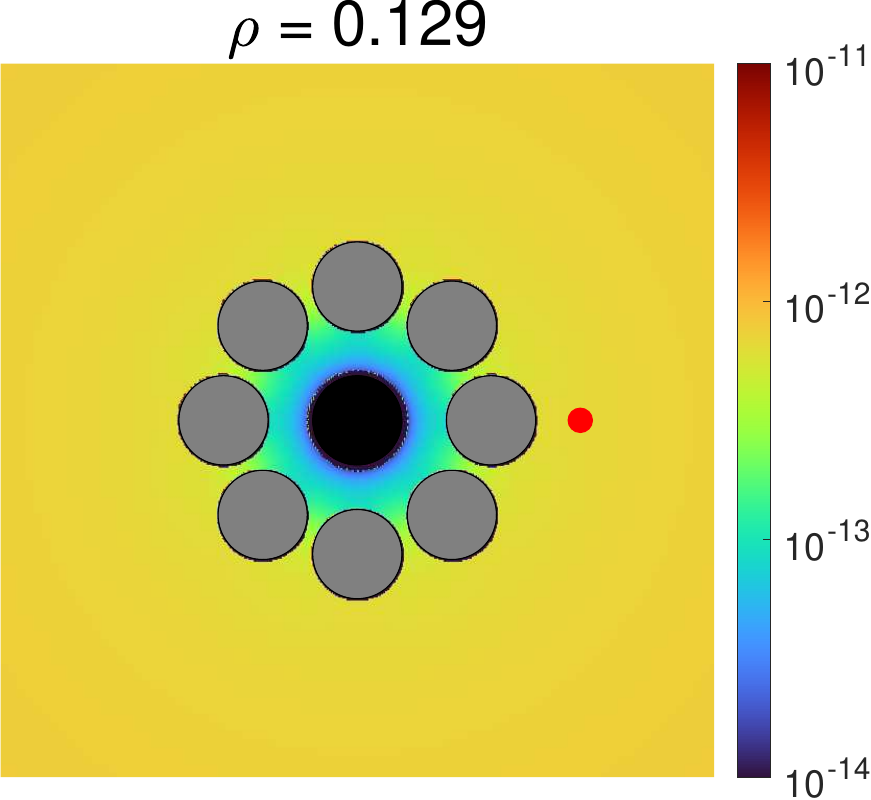}
\includegraphics[height=0.23\textwidth,trim=0pt 0pt 75pt
0pt, clip=true]{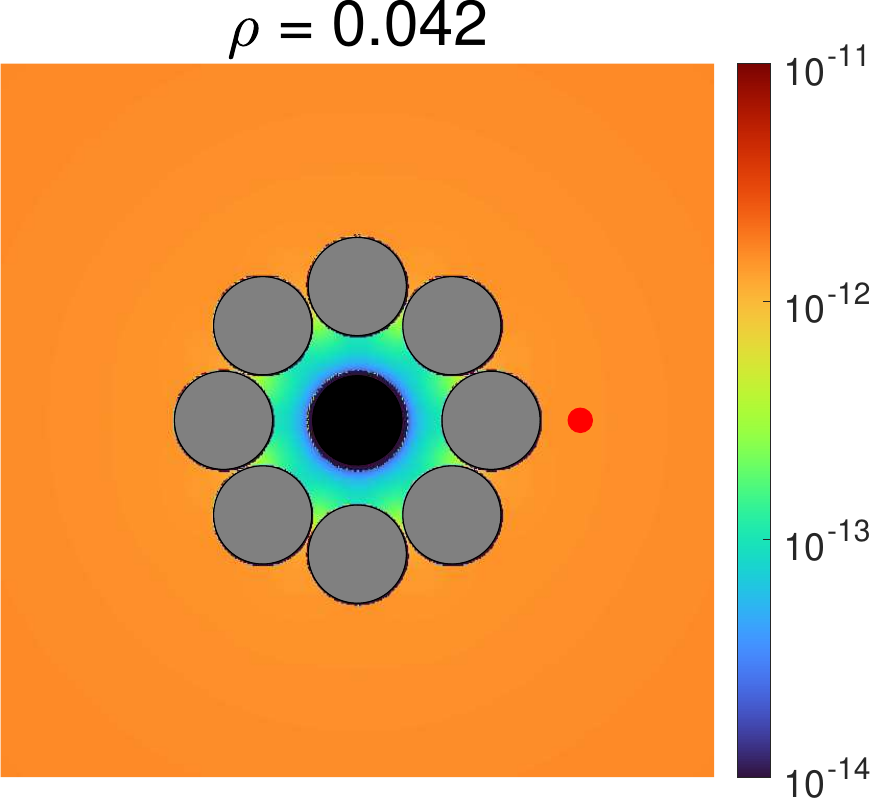}
\includegraphics[height=0.23\textwidth]{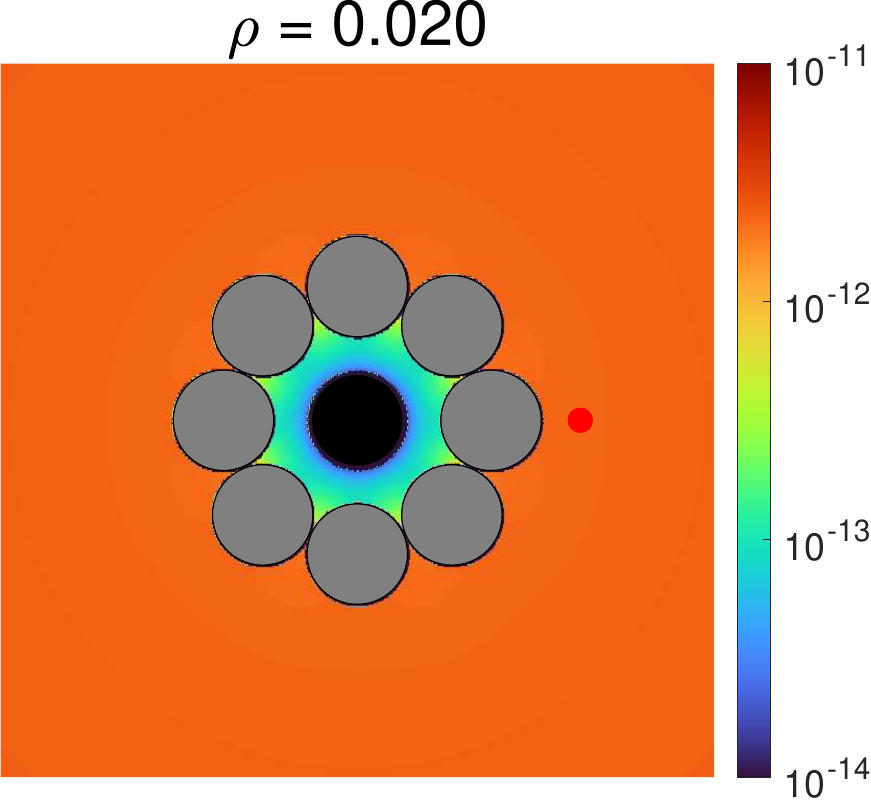}
\caption{\label{fig:faradayHeatLateTime} Heat maps of solutions
$p(\xx,t)$ to \eqref{eqn:diffusion} at $t = 10^{10}$ for geometries with
confining ratios $\rho = \{0.347,0.129,0.042,0.020 \}$ At the highest level of shielding ($\rho = 0.020$), we see
that a particle is roughly $100$ times more likely to be found at a
particular point outside the ring of reflecting bodies of the cage than
at a point inside the ring of reflecting bodies.}
\end{figure}

A heat map of $p(\xx,t)$ is shown in Figure~\ref{fig:faradayHeat} for
each of the four non-zero radii of the reflecting bodies, and at three
values of $t$. Because the Talbot contour scales with $t^{-1}$, the
numerical solution at the earliest time has enough Talbot quadrature
error that a halo appears around the support of the solution (left
column). However, at later times, the numerical solution is very smooth
(right column). At these three early times, the difference between the
solutions is slight, in particular we see small differences in the
density near the gap between reflectors. However, our method allows us
to compute the solution at much later times, a calculation that would be
very challenging using particle based Monte Carlo methods. In
Figure~\ref{fig:faradayHeatLateTime}, we plot the solution at time $t =
10^{10}$. At the highest level of shielding ($\rho = 0.020$), we see
that a particle is about 100 times more likely to be found at a
particular point outside the ring of reflecting bodies of the cage than
at a point inside the ring of reflecting bodies.

\subsection{Multiple Absorbing Traps and Escape from a Maze Geometry}
\label{sec:maze}
\begin{figure}[htp]
\centering
\includegraphics[width=0.5\textwidth]{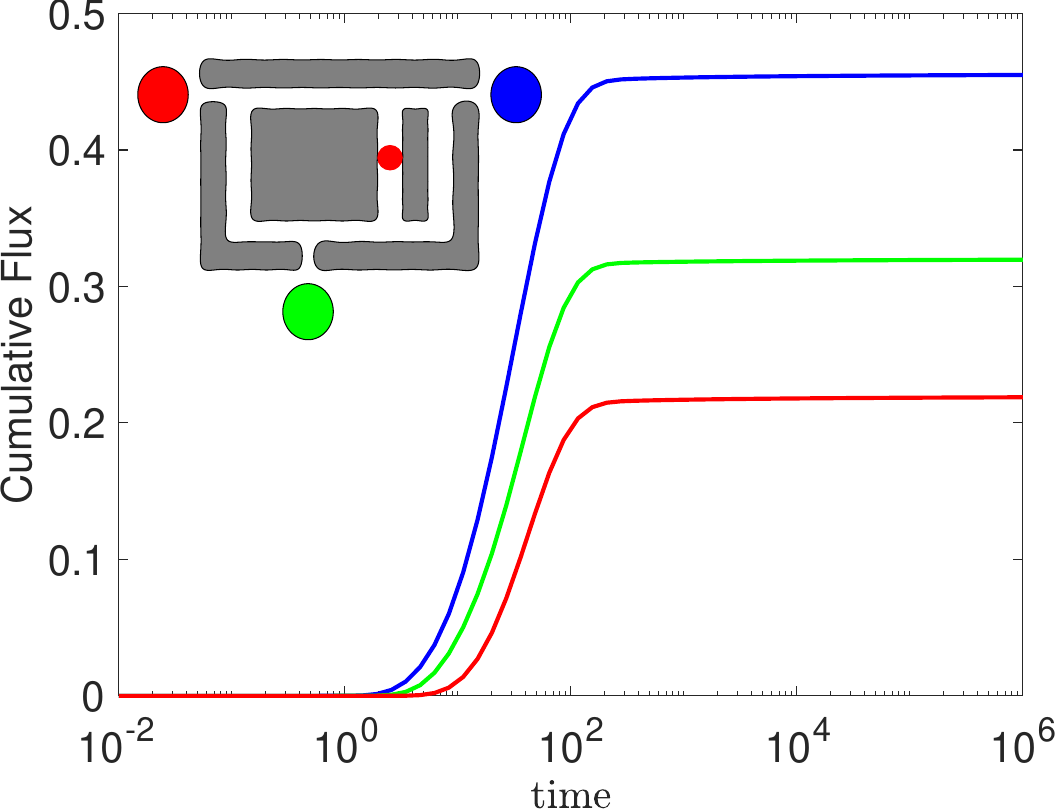}
\caption{\label{fig:mazeCDF} The cumulative flux into three absorbing
traps (colored disks) surrounding a maze of reflecting bodies (gray
bodies) with the initialization point at the center of the maze (red
dot). The colors of the cumulative flux curves match those of the
absorbing disks. Our numerical results show that a diffusing particle
will most likely be absorbed by the blue exit, followed by the green
exit, and is least likely to be absorbed by the red exit.}
\end{figure}

The Faraday cage example involves multiple reflecting bodies, but only a
single absorbing body. In this example, we consider escape from a {\em
maze-like} geometry formed by rectangular blocks with reflecting
boundary conditions. We consider a particle originating at the center of
the maze which diffuses towards one of three exits. The likelihood of
taking a particular exit from the maze is evaluated by calculating the
cumulative flux at absorbing disks placed at each of these apertures. A
schematic of the geometric setup is shown in the inset of
Figure~\ref{fig:mazeCDF}.

\begin{figure}[htp]
\centering
\includegraphics[height=0.25\textwidth,trim=0pt 0pt 75pt
0pt, clip=true]{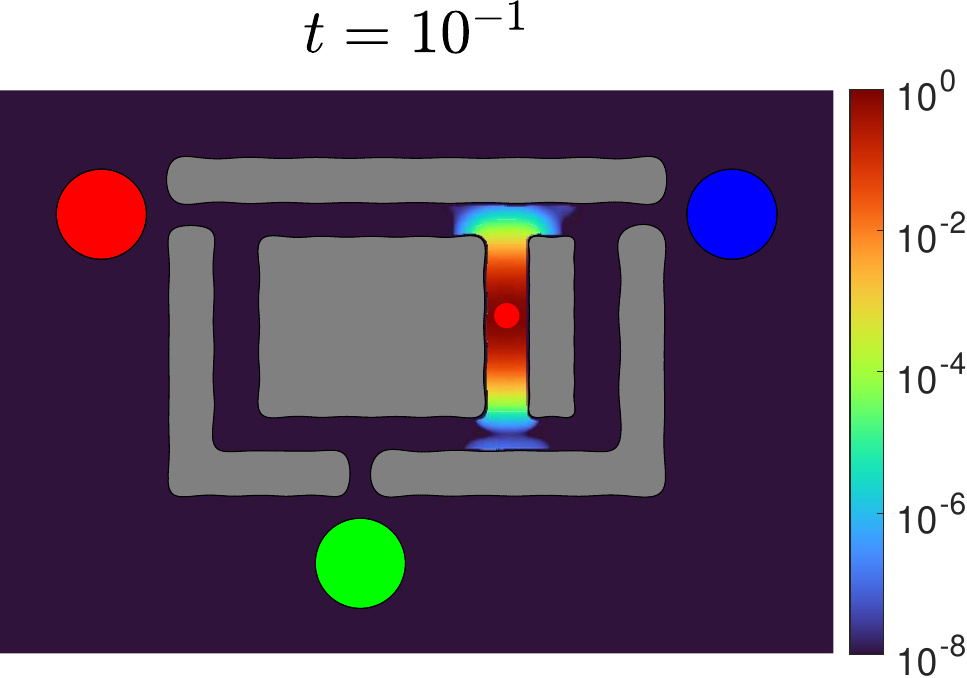}
\includegraphics[height=0.25\textwidth,trim=0pt 0pt 75pt
0pt, clip=true]{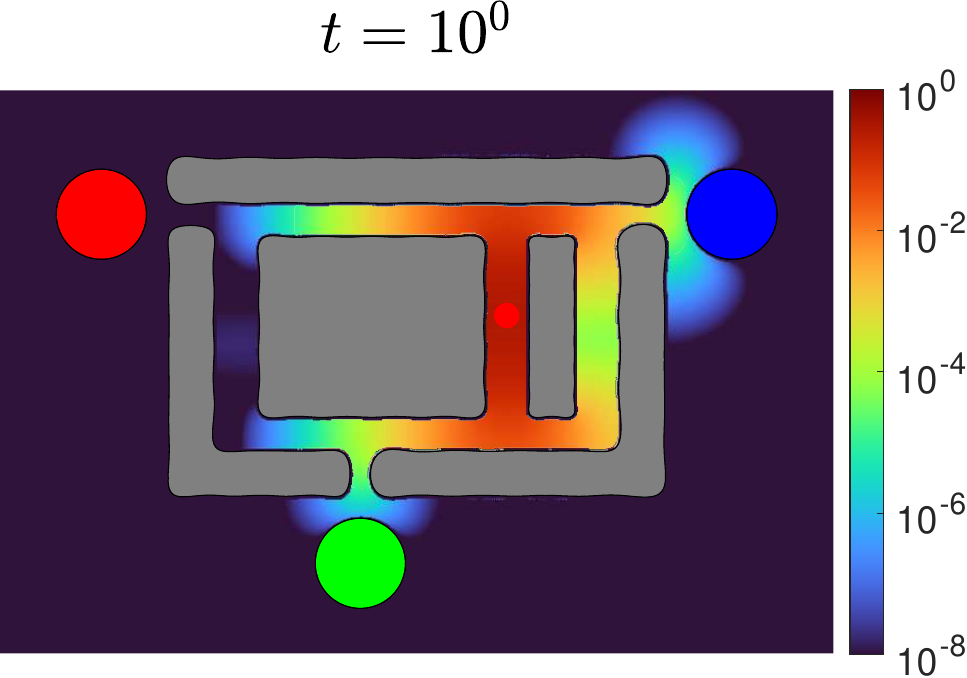}
\includegraphics[height=0.25\textwidth]{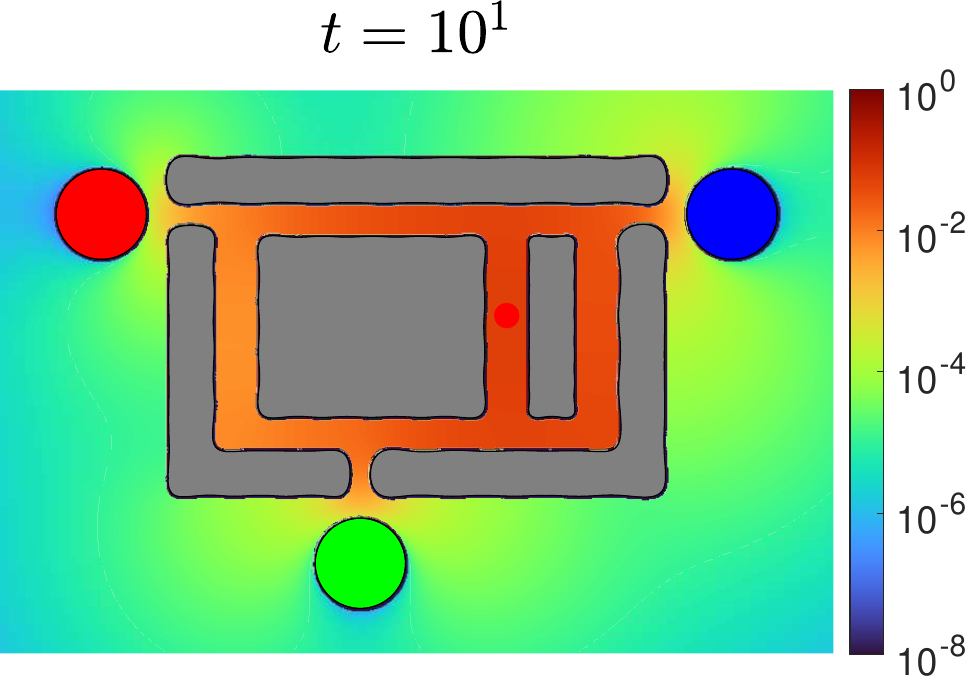}
\caption{\label{fig:mazeHeat} The heat map of solutions $p(\xx,t)$ to
\eqref{eqn:diffusion} at three early times $t =\{10^{-1},10^0,10^1\}$.
The gray polygons are reflecting, the colored disks are absorbers
denoting the maze exits, and the particle is initialized at the maze
center (red dot).}
\end{figure}

In the dynamic results shown in Figure~\ref{fig:mazeCDF}, we see that
for a diffusing particle initialized at the point shown, the blue exit
is the most likely to absorb a particle with probability 45.8\%, the
green exit is the second most likely with probability 32.2\%, and the
red exit is the least likely with probability 22.0\%. These quantities,
which describe the likelihood of hitting any particular absorber, are
known as the \emph{splitting probabilities} and are important in
ecological and signaling applications~\cite{Lindsay2023a, LLM2020,
Venu2015}. In this example, we observe a much faster equilibration
timescale due to the fact that the geometry of the maze funnels most
density close to one of the three absorbers. In terms of particle
trajectories, this means that very few particles escape the maze and the
absorbers while transitioning the exits. Heat maps of $p(\xx,t)$ are
shown in Figure~\ref{fig:mazeHeat}.

\section{Conclusions\label{s:conclusions}}
We have introduced a method for the numerical solution of planar
diffusion problems through a coupling of boundary integral methods with
the Laplace transform. The main application of the approach is to
determine the capture rate of diffusing particles to absorbing sites
through complex geometries. In particular, the methods presented here
can accommodate unbounded geometries comprised of entirely Dirichlet and
entirely Neumann disjoint elements, a scenario that poses a particular
challenge to particle-based methods which can become stuck on reflecting
sections.

Future work should extend this approach to bodies with more general
boundary conditions, including the Robin case and the scenario of a
mixture of Neumann and Dirichlet. The scenario of mixed boundary
conditions arises in cellular signaling where the role of receptors is
modeled through the inclusion of numerous small Dirichlet windows in an
otherwise reflecting body~\cite{BL2018, LWB2017, KAYE2020, LLM2020}. The
case of mixed conditions has been explored previously for the Helmholtz
equation~\cite{Gillman2017} and it would be natural to combine that
approach with the Laplace transform approach of the present work to
enable an extension to parabolic problems. The Robin boundary condition
arises when replacing the mixed boundary conditions by an effective one
of form
\begin{equation}\label{eq:Robin}
  \pderiv{p}{\nn} = \sigma \, p, \qquad \xx \in \Gamma_R.
\end{equation}
The replacement of a mixture of Dirichlet and Neumann boundary
conditions by a uniform Robin condition~\eqref{eq:Robin} is known as
\emph{boundary homogenization}~\cite{Lawley2024, LBS2018, muratov}. This
reduction is a great simplification and is accurate provided the
permeability parameter $\sigma$ is chosen in a way to mimic the capture
rate of the bodies $\Gamma$.

Another important extension of the present work is to handle more
general initial conditions. We remark that our method relies on a known
particular solution of the transform equation to arrive a homogeneous
problem which is solvable by boundary integral means. In the scenario
where the initial condition is a point source, or a summation of point
sources, this particular solution is representable as a combination of
fundamental solutions. For extensions to more general initial
conditions, this particular solution can be represented as a
convolution, however, this introduces an area integral requiring careful
evaluation. Extensions to moving geometries~\cite{Lindsay2017Moving,
TK2015, Lindsay2015} and eroding geometries~\cite{moo-che-chi-qua2023,
qua-moo2018, chi-moo-qua2020}, which can be simulated with BIEs, are
additional important directions for future studies. 

\paragraph{\bf Acknowledgments} A.E.L.~acknowledges support under
National Science Foundation award DMS 2052636. B.D.Q.~acknowledges
support under National Science Foundation award DMS 2012560.

\bibliographystyle{plainnat} 
\biboptions{sort&compress}

\bibliography{refs}

\end{document}